%% file: 0.tex
\documentclass{article}

\usepackage{diagrams}
\usepackage{amsmath}
\usepackage[all,tile]{xy}
\usepackage{hyperref}

\usepackage[english]{mytheorems}
\def\shortversion{}
\usepackage{longversion}
\usepackage{order}
\usepackage{logic}

\input macros.tex

\author{Luigi Santocanale\\
  Laboratoire d'Informatique Fondamentale de Marseille \\
  Universit\'e de Provence \\
  \texttt{luigi.santocanale@lif.univ-mrs.fr}
}
\date{}

\title{Derived Semidistributive Lattices}

\begin{document}

\maketitle

\input abstract.tex
\input introduction.tex
\input general.tex
\input pushdown.tex

\input semidistributive.tex
\input bounded.tex
\input derived.tex
\input examples.tex


\bibliographystyle{elsart-num-sort}
\bibliography{biblio}

\end{document}

%% file: abstract.tex
\begin{abstract}
  For $L$ a finite lattice, let $\Cov(L) \subseteq L^{2}$ denote the
  set of pairs $\gamma = (\gamma_{0},\gamma_{1})$ such that
  $\gamma_{0} \lcover \gamma_{1}$ and order it as follows: $\gamma
  \leq \delta$ iff $\gamma_{0} \leq \delta_{0}$, $\gamma_{1} \not\leq
  \delta_{0}$, and $\gamma_{1} \leq \delta_{1}$.  Let $\Cov(L,\gamma)$
  denote the connected component of $\gamma$ in this poset. Our main
  result states that, for any $\gamma$, $\Cov(L,\gamma)$ is a
  semidistributive lattice if $L$ is semidistributive, and that
  $\Cov(L,\gamma)$ is a bounded lattice if $L$ is bounded.
  Let $\mathcal{S}_{n}$ be the Permutohedron on $n$ letters and
  let $\mathcal{T}_{n}$ be the Associahedron on $n+1$ letters. Explicit
  computations show that $\Cov(\mathcal{S}_{n},\alpha) =
  \mathcal{S}_{n-1}$ and $\Cov(\mathcal{T}_{n},\alpha) =
  \mathcal{T}_{n-1}$, up to isomorphism, whenever $\alpha_{1}$ is an
  atom of ${\cal S}_{n}$ or ${\cal T}_{n}$.

  These results are consequences of new characterizations of finite
  join-semidistributive and of finite lower bounded lattices: (i) a
  finite lattice is join-semidistributive if and only if the
  projection sending $\gamma \in \Cov(L)$ to $\gamma_{0} \in L$
  creates pullbacks, (ii) a finite join-semidistributive lattice is
  lower bounded if and only if it has a strict facet labelling.
  Strict facet labellings, as defined here, are a generalization of
  the tools used by Barbut et al. \cite{BCM} to prove that lattices of
  finite Coxeter groups are bounded.
\end{abstract}


%% file: introduction.tex

\section{Introduction}

The set of covers of a finite lattice comes with a natural ordering
induced by perspectivity.
A cover of a lattice $L$ is an ordered pair $\gamma =
(\gamma_{0},\gamma_{1}) \in L^{2}$ such that the interval
$[\gamma_{0},\gamma_{1}] = \set{x \in L \mid \gamma_{0} \leq x \leq
  \gamma_{1}}$ is the two elements set $\{\gamma_{0},\gamma_{1}\}$. As
usual, we write $\gamma_{0}\lcover \gamma_{1}$ to mean that
$(\gamma_{0},\gamma_{1})$ is a cover.  Two intervals $[x,y]$ and
$[z,w]$ are perspective if either $x = y \land z$ and $w = y \vee z$,
or, vice-versa, $z = x \land w$ and $y = x \vee w$.
We order covers as follows: $\gamma \leq \delta$ if $\gamma_{0} =
\gamma_{1} \land \delta_{0}$ and $\delta_{1} = \gamma_{1} \vee
\delta_{0}$. Thus two covers are comparable if and only if they give
rise to perspective intervals.
 The resulting poset, denoted here by
$\Cov(L)$, is the object investigated in this paper.

\medskip

The main result we shall present is that, whenever $L$ is a
finite
semidistributive lattice, the poset of covers $\Cov(L)$ is the
disjoint union of connected components each of which is again a
semidistributive lattice; if moreover $L$ is a bounded lattice in the
sense of \cite{mckenzie}, then each such component is a bounded
lattice as well. If $\gamma$ is a cover of $L$, then we shall denote
by $\Cov(L,\gamma)$ the connected component of $\gamma$ in $\Cov(L)$
and call it the lattice derived from $L$ by means of 
$\gamma$.  Thus, if $L$ is semidistributive, this process of
constructing derivatives may be iterated.

\medskip

These results are consequences of new characterizations of 
finite
join-semidis\-tributive 
lattices and of finite lower bounded lattices that strengthen well
known facts.  We remark here that, throughout this paper, we shall be
interested in finite lattices only. For this reason and unless
explicitly stated, the word lattice shall be a synonym of finite
lattice.

\smallskip

On one side, it is well known that a lattice is join-semidistributive
if and only if, given a cover $\gamma_{0}\lcover \gamma_{1}$, there
exists a unique cover $j_{\ast}\lcover j$, perspective to
$\gamma_{0}\lcover \gamma_{1}$, such that $j$ is join-irreducible
\cite[\S 2.56]{freese}. The latter property may be rephrased by saying
that for each $\gamma \in \Cov(L)$ there exists a unique $\iota \in
\Cov(L)$, minimal within $\Cov(L)$, such that $\iota$ and $\gamma$ are
comparable.  Following a suggestion of \cite[Theorem 1]{BCM}, we
observe that such uniqueness property is consequence of the pushdown
relation $\pushdown$ between covers being confluent.\footnote{%
  Let us recall that a relation $\rightarrow \subseteq V\times V$ is
  confluent if $v_{0} \rightarrow^{\ast} v_{i}$, $i = 1,2$, implies
  that $v_{i} \rightarrow^{\ast} v_{3}$, $i = 1,2$ for some $v_{3} \in
  V$; here $\rightarrow^{\ast}$ denotes the reflexive transitive
  closure of $\rightarrow$. If the relation $\rightarrow$ has no
  infinite path, then it is a standard result that for each $v$ there
  exists a unique $v'$ such that $v \rightarrow^{\ast} v'$ and $v'$
  has no successor.  } 
This relation shall be introduced in Section \ref{sec:pushdown}; by
studying further it, we refine the existing characterization of
join-semidistributivity to the following
statement: \emph{a lattice is join-semidistributive if
  and only if the poset of covers has pullbacks}.

\smallskip

\begin{figure}[h]
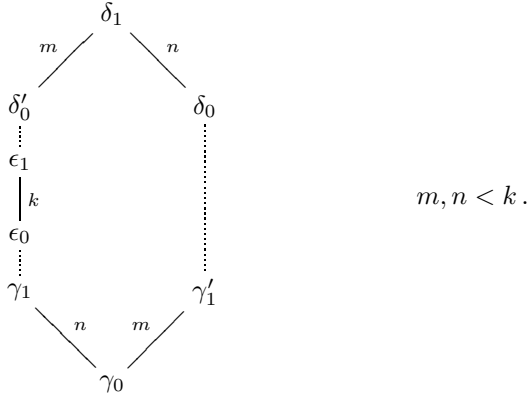

  \centering
  \begin{align*}
    &
    \mydiagram[3.5em]{
      []*+{\delta_{1}}
      (
      -[dl]*+{\delta'_{0}}="A"_{m},
      -[dr]*+{\delta_{0}}="B"^{n}
      ),
      [d(3)]*+{\gamma_{0}}
      (
      -[ul]*+{\gamma_{1}}="C"_{n},
      -[ur]*+{\gamma'_{1}}="D"^{m}
      )
      "B":@{{.}{.}{.}}"D"
      "A":@{{.}{.}{.}}[d(0.6)]*+{\epsilon_{1}}
      -[d(0.8)]*+{\epsilon_{0}}^{k}
      :@{{.}{.}{.}}"C"
    }
    &
    m,n & < k\,.
  \end{align*}
  \caption{A strict labelling of a facet}
  \label{fig:labfacet}
\end{figure}

On the side of lower bounded lattices, we build on the ideas used in
\cite{BCM} to prove that lattices arising from Cayley graphs of finite
Coxeter groups are bounded.  With respect to that work, we move from a
sufficient condition to a complete characterization, and from
boundedness to the weaker notion of lower boundedness. The statement
reads as follows: \emph{a lattice is lower bounded if and only if it
  is join-semidistributive and has a strict facet labelling}.  A
strict facet labelling is a labelling of covers by natural numbers
subject to some constraints. We illustrate next these constraints
under the simplifying assumption that $L$ is a semidistributive
lattice.  As illustrated in Figure \ref{fig:labfacet}, a
facet\footnote{%
  In \cite{BCM} a facet is called a $2$-facet, using a more precise
  wording from combinatorial geometry.  }  is a quadruple of covers
$\delta,\delta',\gamma,\gamma' \in \Cov(L)$ such that $\gamma_{0} =
\gamma'_{0} = \delta_{0} \land \delta'_{0} < \gamma_{1} \vee
\gamma'_{1} = \delta_{1} = \delta'_{1}$, with $\gamma_{1} \leq
\delta'_{0}$ and $\gamma_{1}' \leq \delta_{0}$.  For such a facet, we
shall prove that $\gamma \lcover \delta$ (as well as $\gamma' \lcover
\delta'$) in the poset $\Cov(L)$ and, moreover, that every cover of
covers arises from a facet.  A strict facet labelling assigns the same
number to $\gamma$ and $\delta$ and, consequently, it is constant on
connected components of $\Cov(L)$.  Moreover, such a labelling is
required to be strictly increasing at the interior of a facet which is
a pentagon: if $\epsilon \in \Cov(L)$ is such that $\gamma_{1} \leq
\epsilon_{0} \lcover \epsilon_{1} \leq \delta'_{0}$, then $\delta_{0}$,
$\delta'_{0}$, and $\gamma_{1}$ generate a sublattice which is a
pentagon; then the label of $\epsilon$ should be striclty greater than
the labels of $\gamma$ and of $\delta'$.

\medskip 

We shall carry out some explicit computations: we prove that for the
Permutohedron ${\cal S}_{n}$ -- i.e. the lattice of permutations on
$n$ letters with the weak Bruhat order -- and the Associahedron ${\cal
  T}_{n}$ -- the Tamari lattice on $n+1$ letters -- the relations
\begin{align*}
  \Cov({\cal S}_{n},\alpha) & = {\cal S}_{n -1}\,,
  &
  \Cov({\cal T}_{n},\alpha) & = {\cal T}_{n -1}\,,
\end{align*}
hold up to isomorphism, whenever $\alpha$ is of the form
$(\bot,\alpha_{1})$, that is, $\alpha_{1}$ is an atom.  A similar
relation holds for ${\cal B}^{n}$, the Boolean algebra on $n$
atoms. Since in ${\cal B}^{n}$ the join irreducible elements are
exactly the atoms, the property is there stronger in that
$\Cov({\cal B}^{n},\gamma) = {\cal B}^{n-1}$ \emph{for every $\gamma
  \in \Cov({\cal B}^{n})$}. As a consequence 
\begin{align}
  \label{eq:nBn-1}
  \Cov({\cal B}^{n}) = \coprod_{i=1,\ldots ,n} {\cal B}^{n-1} = n
  \cdot {\cal B}^{n-1}\,,
\end{align}
where the coproduct and the products -- implicit in the exponents $i$
of ${\cal B}^{i}$ -- are taken within the category of posets and order
preserving functions.  A reason for the name of derived lattices
becomes transparent if, in equation \eqref{eq:nBn-1}, we replace the
symbol $\Cov$ with the symbol $\frac{\partial}{\partial {\cal B}}$.

\medskip

It is suggestive to call the lattices ${\cal S}_{n}, {\cal T}_{n}$,
and ${\cal B}^{n}$ \emph{regular}, meaning that the shape of $
\Cov(L,\alpha)$ does not depend on the choice of the atom
$\alpha_{1}$.  This use of terminology from combinatorial geometry is
on purpose since, with this work, we aim at giving a ground to some
intuitions relating algebra, order, and geometry. In combinatorial
geometry a Permutohedron is a particular convex polytope
\cite[0.10]{ziegler}. By orienting the graph reduct of this polytope
we obtain the Hasse diagram of the lattice of permutations.
We shall show that the 2-facets of this polytope can also be oriented
so that they give rise to the Hasse diagram of $\Cov({\cal S}_{n})$,
i.e. to covers of covers, or say $2$-covers.
By semidistributivity, we can define $n$-covers as elements of
$\Cov^{n}({\cal S}_{n})$. It is not difficult to verify that these
correspond to oriented $n$-facets of the polytope ${\cal S}_{n}$.

Since we can define $n$-facets for an arbitrary semidistributive
lattice, it becomes natural to ask whether the lattice theoretic
algebra plays any role in geometry and, for example, whether we can
characterize lattices that arise from convex polytopes such as ${\cal
  S}_{n},{\cal T}_{n}$, and ${\cal B}^{n}$ by lattice theoretic means.
This paper will not answer these questions, its purpose being limited
to settling a ground for future researches.

\medskip

Let us also stress on the fact that ideas and results presented here
have their origin at the intersection between order theory and the
theory of rewriting systems, in the spirit of \cite{newman} and
\cite{bb}.  Even if we are not going to emphasize this aspect, it is
worth recalling it.  It was suggested in \cite{BCM} that the
perspectivity relation between covers is generated by a sort of
rewriting system. In section \ref{sec:pushdown} we make explicit such
a rewriting system, namely, we define the pushdown relation
$\pushdown$ between covers. We also explicitly introduce the class of
pushdown lattices, as our proofs make heavy use its properties.
Briefly, a lattice is pushdown if the perspectivity order on covers is
generated by the pushdown relation, thus allowing a sort of local
reasoning on the global structure of the covers.

\medskip

The paper is structured as follows. We first recall the basic concepts
concerning lattices in Section \ref{sec:general} and we introduce
pushdown lattices in Section \ref{sec:pushdown}. We give in Section
\ref{sec:jsem} our characterization of join-semidistributive lattices,
and in Section \ref{sec:bounded} our characterization of lower bounded
lattices. We define in Section \ref{sec:derived} lattices of the form
$\Cov(L,\gamma)$, that is lattices derived from semidistributive
lattices, and prove then that properties such as semidistributivity
and being lower bounded lift from $L$ to $\Cov(L,\gamma)$. Finally, in
Section \ref{sec:examples}, we exemplify the construction of derived
lattices on Newman lattices.


%% file: general.tex
\section{Preliminaries}
\label{sec:preliminaries} 
\label{sec:general} 

We begin by introducing standard definitions and notations on
finite lattices.

\smallskip

Let $P$ be a poset.  A \emph{cover} in $P$ is an ordered pair
$(\gamma_{0},\gamma_{1}) \in P^{2}$ such that 
$\gamma_{0} \leq \gamma_{1}$ and the closed interval $\set{x \in
  L\mid\gamma_{0} \leq x \leq \gamma_{1}}$ is the two element set
$\set{\gamma_{0},\gamma_{1}}$ -- in particular $\gamma_{0} \neq
\gamma_{1}$.  As usual, we shall write $\gamma_{0} \lcover \gamma_{1}$
if $(\gamma_{0},\gamma_{1})$ is a cover and say that $\gamma_{0}$ is a
lower cover of $\gamma_{1}$ and that $\gamma_{1}$ is an upper cover of
$\gamma_{0}$. We shall denote by $\Cov(P)$ the set of covers of $P$
and use Greek letters $\gamma,\delta \ldots $ to denote these covers.

\smallskip

An element of a lattice $L$ is \emph{join-}(resp.
\emph{meet-})\emph{irreducible} if it has a unique lower (resp.
upper) cover.  We denote by $J(L)$ (resp. $M(L)$) the set of
join-(resp. meet-)irreducible elements of $L$. If $j \in J(L)$ then
$j_{\ast}$ denotes the unique element of $L$ such that $j_{\ast}
\lcover j$.  If $m \in M(L)$, then $m^{\ast}$ denotes the unique
element of $L$ such that $m \lcover m^{\ast}$. Let us introduce the
standard arrow relations between join-irreducible and meet-irreducible
elements.  For $j \in J(L)$ and $m \in M(L)$, we write $j \jmup m$ if
$j \leq m^{\ast}$ and $j \not \leq m$, and $j \jmdown m$ if $j_{\ast}
\leq m$ and $j \not\leq m$; we write $j \jmpersp m$
if $j \jmup m$ and $j \jmdown m$.

\smallskip

We finally recall that a lattice is \emph{join-semidistributive} if it
satisfies the Horn sentence
\begin{align}
  \label{eq:sdvee}
  \tag{$SD_{\vee}$}
  x \vee y  = x \vee z & \Rightarrow x \vee (y \land z) = x\vee y\,.
\end{align}
A lattice is \emph{meet-semidistributive} if it satisfies the Horn
sentence dual of \eqref{eq:sdvee}. It is \emph{semidistributive} if
it is both meet-semidistributive and join-semidistributive.

\paragraph{Posets with pullbacks.}
\label{sec:prelpbs} 
These posets will play a central role in our development. We introduce
them now together with their elementary properties.  Recall from
\cite{BCM} that a \emph{hat} in a poset $P$ is a triple $(u,v,w)$ such
that $u \lcover v$, $w \lcover v$, and $u\neq
w$. 
An \emph{antihat} in $P$ is defined dually.
A \emph{cospan} in $P$ is a triple of elements $(u,v,w)$ such that $u
\leq v$ and $w \leq v$; in particular a hat is a particular kind of a
cospan. We say that a cospan $(u,v,w)$ \emph{has a pullback} if the
meet $u \land w$ exists.
\begin{definition}
  \label{def:pospbs}
  We say that a poset \emph{$P$ has pullbacks} if every cospan
  in $P$ has a pullback. We say that $f: P \rTo Q$ \emph{preserves
    pullbacks} if $f(u \land w)$ is the pullback of the cospan
  $(f(u),f(v),f(w))$, whenever $u \land w$ is the pullback of the
  cospan $(u,v,w)$.
\end{definition}
Clearly, a poset has pullbacks iff every finite non empty set
admitting an upper bound has a meet.  Notice that every
meet-semilattice is a poset with pullbacks. The following diagram
exhibits a poset with pullbacks which is not a meet-semilattice.
$$
\mydiagram[3em]{
  []*+{\bullet}="A"
  (-[d]*+{\bullet}="C"
  -[dr]*+{\bullet}="E"
  )
  [rr]*+{\bullet}="B"
  -[d]*+{\bullet}="D"
  -"E"
  "A"-"D"
  "B"-"C"
}
$$
The following Proposition is almost a reformulation of Definition
\ref{def:pospbs}.
\begin{proposition}
  A poset $P$ has pullbacks iff every principal ideal of $P$ is a
  lattice.
\end{proposition}
We state next, without proofs, some facts illustrating the specific
role of pullbacks of hats among all the pullbacks.
\begin{lemma}
  \label{lemma:pullbacksofhats}
  A finite poset $P$ has pullbacks iff every hat has a
  pullback.
\end{lemma}
\begin{longversion}
  \begin{proof}
    One direction is obvious. Let us assume that every hat has a
    pullback.  The height of a cospan $(u,v,z)$ is the height of $v$.
    We prove that every cospan has a pullback by induction on the
    height.

    Let us consider the cospan $(u,v,z)$ and let us assume that every
    cospan of height less than $h(v)$ has a pullback.

    If $u = v$ or $z = v$ then the meet of $u$ and $z$ exists, since
    $u$ and $z$ are comparable.

    Otherwise, let $u \leq u' \lcover v$ and $z \leq z' \lcover v$.
    If $u' = z'$, then $u \land z$ exists since the height of the
    cospan $(u,u',z)$ is less than $h(v)$. Otherwise $u'\neq z'$ so
    that $(u',v,z')$ is a hat and $u ' \land z'$ exists.  By the
    inductive hypothesis $u \land u' \land z'$ and $z \land u' \land
    z'$ exists.  Also, the cospan $(u \land u' \land z',u' \land z',z
    \land u' \land z')$ has a pullback $z \land u \land u' \land z'$,
    which by routine verification it is shown to be the meet $u \land
    z$.
  \end{proof}
\end{longversion}
\begin{lemma}
  Let $P,Q$ be finite posets with pullbacks. If $f: P \rTo Q$
  preserves pullbacks of hats, then it preserves pullbacks.
\end{lemma}
\begin{longversion}
  \begin{proof}
    We prove that the pullback of a cospan is preserved by $f$ by
    induction on the height.

    Let us consider the cospan $(u,v,z)$ and let us assume that the
    pullback of every cospan of height less than $h(v)$ is preserved.

    Clearly this is the case if $u = v$ or $z = v$, since then $f(u)$
    and $f(z)$ are comparable.

    Otherwise, let $u \leq u' \lcover v$ and $z \leq z' \lcover v$.
    If $u' = z'$, then $f(u \land z) = f(u) \land f(z)$ since the
    height of the cospan $(u,u',z)$ is less than $h(v)$. Otherwise
    $u'\neq z'$ so that $(u',v,z')$ is a hat and $f(u ' \land z') =
    f(u') \land f(z')$.  By the inductive hypothesis $f(u \land u'
    \land z') = f(u) \land f(u' \land z')$ and $f(z \land u' \land z')
    = f(u) \land f(u' \land z')$ and $f(z \land u' \land z') = f(u)
    \land f(u' \land z')$ and $f(z) \land f(u' \land z')$.  The cospan
    $(u \land u' \land z',u' \land z',z \land u' \land z')$ has height
    less than $h(v)$, hence its pullback $z \land u \land u' \land z'$
    is preserved by $f$.  Therefore we get
    \begin{align*}
      f(u \land z) & = f(z \land u \land u' \land z',u' \land z')
      \tag*{by the previous
        Proposition}\\
      & = f(u \land u' \land z') \land f(z \land u' \land z') \\
      & = f(u) \land f(u' \land z') \land f(z) \\
      & = f(u) \land f(u') \land f(z') \land f(z) \\
      & = f(u) \land f(z)\,.
    \end{align*}
  \end{proof}
\end{longversion}

The following observations will turn to be more interesting for our
goals.
\begin{proposition}
  If a finite poset $P$ has pullbacks, then each connected component
  of $P$ has a least element.
\end{proposition}
\begin{proof}
  Since a principal ideal is a finite lattice, then it has a least
  element. Let therefore $\bot_{x}$ denote the least element of the
  principal ideal of $x$.  We argue that if $x,y$ belong to the same
  connected component, then $\bot_{x} = \bot_{y}$.  To this goal it
  suffices to establish that $\bot_{x} = \bot_{y}$ whenever $z \leq
  x,y$ for some $z$. We have $\bot_{x} \leq z \leq y$, and hence
  $\bot_{y} \leq \bot_{x}$. By symmetry, $\bot_{x} \leq \bot_{y}$.
\end{proof}

Let us say that $P$ \emph{has pushouts} if its dual poset $P^{op}$ has
pullbacks.
\begin{corollary}
  \label{cor:pullpush}
  If a finite poset $P$ has pullbacks and pushouts, then each
  connected component of $P$ is a lattice.
\end{corollary}
\begin{proof}
  Since $P$ has pushouts each connected component has a maximum
  element, that is, it is a principal ideal. Since $P$ has pullbacks,
  such an ideal is a lattice.
\end{proof}


%% file: pushdown.tex
\section{Pushdown Lattices}
\label{sec:pushdown}

In the following $L$ will denote a fixed 
lattice.
Let us recall that an interval $[x,y]$ is said to \emph{transpose
  down} to an interval $[z,w]$ if $z = x \land w$ and $y = x \vee
w$. Two intervals $[x,y]$ and $[z,w]$ are said to be
\emph{perspective} if either $[x,y]$ transposes down to $[z,w]$, or
$[z,w]$ transposes down to $[x,y]$. %
Perspectivity suggests how to define an ordering
$\leq$ on the set $\Cov(L)$ of covers
of $L$: for $\gamma, \delta \in \Cov(L)$, we let
\begin{align*}
  \gamma \leq \delta & \tif
  \gamma_{0} \leq \delta_{0}, \,\gamma_{1} \not\leq \delta_{0}, \tand
  \gamma_{1} \leq \delta_{1}\,.
\end{align*}
It is easy to verify that, for
$\gamma,\delta \in \Cov(L)$, $\gamma \leq \delta$ iff
$[\gamma_{0},\gamma_{1}]$ and $[\delta_{0},\delta_{1}]$ are
perspective with $[\delta_{0},\delta_{1}]$ transposing down to
$[\gamma_{0},\gamma_{1}]$.
This is an order relation on $\Cov(L)$, since $L$ is a lattice: if
$\gamma \leq \delta \leq \epsilon$ then clearly $\gamma_{i} \leq
\epsilon_{i}$ for $i = 0,1$, and if $\gamma_{1} \leq \epsilon_{0}$,
then also $\gamma_{1} \leq \delta_{1} \land \epsilon_{0} =
\delta_{0}$, a contradiction.  Observe also that a sufficient
condition for such a relation to be an ordering is that $L$ has
pullbacks. When referring to the poset $\Cov(L)$ we shall mean the
ordered pair $\langle \Cov(L), \leq \rangle$.  By $\Cov(L,\gamma)$ we
shall denote the connected component of $\gamma$ in $\Cov(L)$.  

\breath
The
two projections
\begin{align*}
  \pr[i] & : \Cov(L) \rTo L\,, & i = 0,1\,,
\end{align*}
sending $\gamma$ to $\gamma_{i}$,  are order preserving. They
will play a key role in the rest of the paper.

Let $P, Q$ be posets, an order preserving map $\pi : P \rTo Q$ is
\emph{conservative} if it strictly preserves the order, i.e. $x \leq
y$ and $\pi(x) = \pi(y)$ imply $x = y$, or, equivalently $x < y$
implies $\pi(x) < \pi(y)$. Our first remark is the following:
\begin{lemma}
  \label{lemma:conservative}
  In any lattice $L$ the projections $\pr[i]$, $i = 0,1$, are
  conservative.
\end{lemma}
\begin{proof}
  If $\gamma \leq \delta$ and $\gamma_{0} = \delta_{0}$, then the
  relations $\delta_{0} = \gamma_{0} \lcover \gamma_{1} \leq
  \delta_{1}$ and $\delta_{0} \lcover \delta_{1}$ imply $\delta_{1} =
  \gamma_{1}$.
\end{proof}

If $P, Q$ are two posets, then an order preserving function $\pi : Q
\rTo P$ is a \emph{Grothendieck fibration} if for each $\delta \in Q$
the restriction of $\pi$ to the principal ideal generated by $\delta$
is an embedding.  Spelled out, this means that $\gamma,\epsilon\leq
\delta$ and $\pi(\gamma) \leq \pi(\epsilon)$ implies $\gamma \leq
\epsilon$.
\begin{definition}
  We say that a lattice $L$ is a \emph{pushdown} lattice if the
  projection $\pr[0] : \Cov(L) \rTo L$ is a Grothendieck fibration.
\end{definition}
We shall say that $L$ is a \emph{pushup} lattice if the dual $L^{op}$
is pushdown.  Spelled out,  $L$ is pushdown if the following property
holds:
\begin{align*}
  \gamma,\epsilon\leq \delta \tand \gamma_{0} \leq \epsilon_{0}
  &
  \text{ imply } \gamma \leq \epsilon \,.
\end{align*}
Later, we shall refer to this property as to the \emph{pushdown
  property}.  Dually, $L$ is pushup if $\delta\leq \gamma,\epsilon$
and $\epsilon_{1} \leq \gamma_{1}$ imply $\epsilon \leq \gamma$. These
properties may be simplified even more, for example $L$ is pushdown
iff $\gamma,\epsilon\leq \delta$ and $\gamma_{0} \leq \epsilon_{0}$
imply $\gamma_{1} \leq \epsilon_{1}$.

\breath

Examples of lattices enjoying these properties arise from some form of
semidistributivity.
\begin{lemma}
  \label{lemma:joinsempushdown}
  If $L$ is a join-semidistributive lattice, then it is both a pushdown
  and a pushup lattice.
\end{lemma}
\begin{proof}
  The proof is sketched in the two diagrams below:
  $$
  \xygraph{
    []!~:{@{.}}
    []!c{\gamma}
    [ru]!c{\epsilon}
    [ru]!c{\delta}
    "\gamma_{0}":"\epsilon_{0}":"\delta_{0}"
    "\epsilon_{1}":"\delta_{1}"
    "\gamma_{1}":@/^2em/"\delta_{1}"
    "\delta_{0}"[r(0.4)]*+{=x}
    "\epsilon_{1}"[r(0.4)]*+{=y}
    "\gamma_{1}"[r(0.4)]*+{=z}
  }
  \hspace{10mm}
  \xygraph{%
    []!~:{@{.}}
    [uu]!c{\gamma}
    [rd]!c{\epsilon}
    [rd]!c{\delta}
    "\gamma_{1}":"\epsilon_{1}":"\delta_{1}"
    "\epsilon_{0}":"\delta_{0}"
    "\gamma_{0}":@/_2em/"\delta_{0}"
    "\gamma_{0}"[r(0.4)]*+{=x}
    "\epsilon_{0}"[r(0.4)]*+{=y}
    "\delta_{1}"[r(0.4)]*+{=z}
  }
  $$
  Let $x,y,z$ be as above on the left: if $\gamma_{1} \not\leq
  \epsilon_{1}$, then $y \land z = \epsilon_{1} \land \gamma_{1} =
  \gamma_{0}$ and, consequently, $x \vee y = x \vee z$ but $x \vee
  (y\land z) = x \vee \gamma_{0} = x < \delta_{1} = x \vee y$.

  Let $x,y,z$ be as above on the right: if $\gamma_{0} \not\geq
  \epsilon_{0}$, then $x \vee y = \gamma_{0} \vee \epsilon_{0} =
  \gamma_{1}$ and, consequently, $x \vee y = x \vee z$ but $x \vee
  (y\land z) = x \vee \delta_{0} = x < \gamma_{1} = x \vee y$.
\end{proof}
By duality, it follows that meet-semidistributive lattices are both
pushdown and pushup lattices.

\breath

To understand the relationship between pushdown lattices and hats and
antihats as defined in Section \ref{sec:prelpbs}, let us introduce the
\emph{pushdown relation $\pushdown$} on $\Cov(L)$ as follows.  For $u \in L$ and
$\gamma,\delta \in \Cov(L)$, let us write $\delta \pushdown[u] \gamma$
if $u \neq \delta_{0}$, $\gamma_{0} = u \land \delta_{0}$, and
$\gamma_{1} \leq u \lcover \delta_{1}$.  For $\gamma,\delta \in
\Cov(L)$, let us write $\delta \pushdown \gamma$ if $\delta
\pushdown[u] \gamma$ for some $u \in L$.
A first remark is that $\delta \pushdown \gamma$ implies $\gamma <
\delta$, in any lattice: if $\delta \pushdown[u] \gamma$, then
$\gamma_{0} = u \land \delta_{0} \leq \delta_{0}$, $\gamma_{1} \leq u
\leq \delta_{1}$; if moreover $\gamma_{1} \leq \delta_{0}$, then
$\gamma_{0} < \gamma_{1} \leq u \land \delta_{0} = \gamma_{0}$, a
contradiction. 
Next, we observe that if $\delta \pushdown[u]\gamma$, then
$(u,\delta_{1},\delta_{0})$ is a hat and $\gamma_{0}$ is its pullback.
The next Lemma implies that, in a pushdown lattice, to a given hat
$(u,\delta_{1},\delta_{0})$ there corresponds a unique antihat
$(x,\gamma_{0},y)$ such that $\gamma_{0} = u \land \delta_{0}$, $x
\leq u$, and $y \leq \delta_{0}$.
\begin{lemma}
  \label{lemma:uniquepusher}
  In a pushdown lattice, for each $\delta \in \Cov(L)$ and $u\in L $
  such that $u \lcover \delta_{1}$ and $u \neq \delta_{0}$, there
  exists a unique $\gamma$ such that $\delta \pushdown[u] \gamma$.
\end{lemma}
\begin{proof}
  Let $u \in L$ and $\delta \in \Cov(L)$ be as in the statement of the
  Lemma.

  We first construct $\gamma \in \Cov(L)$ such that $\delta
  \pushdown[u] \gamma$.  Let $\gamma_{0} = u \land \delta_{0}$ and
  observe that $\gamma_{0} < u$, since $u,\delta_{0}$ are not
  comparable. Hence, we can choose $\gamma_{1}$ such that $\gamma_{0}
  \lcover \gamma_{1} \leq u$ and define $\gamma =
  (\gamma_{0},\gamma_{1})$. It is easily verified that $\delta
  \pushdown[u] \gamma$.
 
  Let us suppose next that $\delta\pushdown[u] \gamma$ and
  $\delta\pushdown[u] \epsilon$, we shall show that $\gamma =
  \epsilon$. By the definition of $\pushdown$, $\gamma_{0}= u \land
  \delta_{0} = \epsilon_{0}$, hence $\gamma_{0} \leq
  \epsilon_{0}$. Also $\delta \pushdown[u] \gamma$ implies $\gamma
  \leq \delta$ and, similarly, $\delta \pushdown[u] \epsilon$ implies
  $\epsilon \leq \delta$. Hence $\gamma,\epsilon \leq \delta$,
  $\gamma_{0} \leq \epsilon_{0}$, and the pushdown property imply
  $\gamma \leq \epsilon$. By symmetry, we obtain $\epsilon \leq
  \gamma$ as well.
\end{proof}

\begin{proposition}
  \label{prop:transcolos}
  In a pushdown lattice $L$ the order of $\Cov(L)$ is the converse
  relation of the reflexive transitive closure of $\pushdown$.
\end{proposition}
\begin{proof}
  Let $\pushdown[*]$ denote the reflexive transitive closure of
  $\pushdown$. We have already seen that $\delta \pushdown \gamma$
  implies $\gamma < \delta$ in any lattice, hence, if $\delta
  \pushdown[*] \gamma$, then $\gamma \leq \delta$ as well. 
  Therefore we shall focus on proving the converse implication, that
  is, if $L$ is pushdown and $\gamma \leq \delta$, then we can find a
  path of the relation $\pushdown$ from $\delta$ to $\gamma$.  The
  proof is by induction on the height of the interval
  $[\gamma_{0},\delta_{0}]$. 

  If $\gamma_{0} = \delta_{0}$, then $\gamma \leq \delta$ and Lemma
  \ref{lemma:conservative} imply $\gamma = \delta$. Thus, the empty
  path witnesses the relation $\delta \pushdown[*] \gamma$.

  Let us suppose that $\gamma_{0} < \delta_{0}$ so that, by Lemma
  \ref{lemma:conservative}, $\gamma_{1} < \delta_{1}$ as well. Pick $u
  \in L$ such that $\gamma_{1} \leq u \lcover \delta_{1}$ and observe
  that $u \neq \delta_{0}$, since $\gamma_{1} \not\leq
  \delta_{0}$. Let $\epsilon$ be determined by the property that
  $\delta \pushdown[u] \epsilon$. Then $\gamma_{0} \leq u \land
  \delta_{0} = \epsilon_{0}$, $\gamma,\epsilon \leq \delta$, and
  therefore, by the pushdown property, $\gamma \leq
  \epsilon$. Moreover, $[\gamma_{0},\epsilon_{0}] \subset
  [\gamma_{0},\delta_{0}]$ and, by the induction hypothesis, $\epsilon
  \pushdown[*] \gamma$. Considering that $\delta \pushdown \epsilon
  \pushdown[*] \gamma$, it follows that $\delta \pushdown[*] \gamma$.
\end{proof}

We end this section with a characterization of pushdown lattices in
terms of the relation $\pushdown$.  This was actually our original
definition of pushdown lattices: if $\delta \pushdown[u] \gamma$,
then the cover $\delta$ has been lowered to the cover $\gamma$ by
means of the push of $u$. This intuition is already present in
\cite{BCM}. %
\begin{proposition}
  \label{prop:equivpushdown}
  A lattice $L$ is a pushdown lattice if and only if $\gamma \leq
  \delta$ and $\delta \pushdown[u] \epsilon$ with $\gamma_{0} \leq u$
  imply $\gamma \leq \epsilon$.
\end{proposition}
\begin{proof}
  Clearly, if a lattice is pushdown, $\gamma \leq \delta$, $\delta
  \pushdown[u] \epsilon$, and $\gamma_{0} \leq u$, then $\gamma_{0}
  \leq u \land \delta_{0} = \epsilon_{0}$, $\gamma,\epsilon \leq
  \delta$, and hence $\gamma \leq \epsilon$.

  Conversely, let us suppose that $\gamma,\epsilon \leq \delta$ with
  $\gamma_{0} \leq \epsilon_{0}$ and that $L$ has the property stated
  in the Proposition. 
  The reader will have no difficulties to adapt the proof of
  Proposition \ref{prop:transcolos} in order to construct a path from
  $\delta$ to $\epsilon$ of the form $\delta =
  \theta^{0}\spushdown[3pt]{u_{1}} \theta^{1}\ldots \theta^{n-1}
  \spushdown[3pt]{u_{n}} \theta^{n} =\epsilon$.
  By assumption we have $\gamma \leq \delta = \theta^{0}$.  Let us
  suppose next that $\gamma \leq \theta^{i}$, for some $i \in
  \set{0,\ldots ,n-1}$.  Since $\theta^{i} \spushdown[3pt]{u_{i +1}}
  \theta^{i +1}$ and $\gamma_{0} \leq \epsilon_{0} \leq u_{i +1}$, the
  property of the Proposition ensures that $\gamma\leq \theta^{i +1}$.
  We have therefore $\gamma \leq \theta^{n} = \epsilon$, showing that
  $L$ is pushdown.
\end{proof}


%% file: semidistributive.tex
\section{Join-Semidistributive Lattices}
\label{sec:jsem}

Many characterizations of (finite) join-semidistributive lattices are
already available, see \cite[Theorem 1.11]{adaricheva} for an example.
In this section we introduce one more characterization, Theorem
\ref{theo:charjsemidis}. A closer glance to this Theorem shows that it
is a refinement of a well known characterization 
stating that a lattice is join-semidistributive if and only if for
each meet-irreducible element $m$ there exists a unique
join-irreducible element $j$ such that $j \jmpersp m$.  This
characterization may be rephrased in terms of the poset $\Cov(L)$. The
correspondence sending a meet-irreducible element $m$ to the cover
$(m,m^{\ast})$ establishes a bijection between $M(L)$ and the set of
maximal elements of $\Cov(L)$. A similar bijection may be defined
between $J(L)$ and the set of minimal elements of $\Cov(L)$. With
these bijections at hand, we can state the previous characterization
as follows:
\begin{lemma}[See Theorem 2.56 in \cite{freese}]
  \label{lemma:usual}
  A lattice $L$ is join-semidistributive if and only if each maximal
  element of $\Cov(L)$ has a least element below it.
\end{lemma}
Roughly
speaking, the new characterizations is obtained from the previous one
by replacing the statement \emph{each maximal element of $\Cov(L)$ has
  a least element below it} with the more informative statement
\emph{$\Cov(L)$ has pullbacks}.




\breath

An order preserving function $\pi: P \rTo Q$ \emph{creates pullbacks}
if the following condition holds: whenever $x,y,z \in P$ are such that
$x ,y \leq z$ and $\pi(x) \land \pi(y)$ exists in $Q$, then there
exists a unique $u \leq x,y$ such that $\pi(u) = \pi(x) \land \pi(y)$;
moreover $u = x \land y$.  For conservative order preserving maps this
condition splits as the conjunction of two conditions, as stated in
the following Lemma.
\begin{lemma}
  \label{lemma:picreatespbs}
  A conservative order preserving function $\pi: P \rTo Q$
  \emph{creates pullbacks} if and only if (i) it is a Grothendieck
  fibration and (ii) if $x ,y \leq z \in P$ and $\pi(x) \land \pi(y)$
  exists in $Q$, then there exists $u \leq x,y$ such that $\pi(u) =
  \pi(x) \land \pi(y)$.
\end{lemma}
\begin{proof}
  Let us suppose that $\pi$ creates pullbacks, so that (ii) 
  certainly holds. If $x,y \leq z$ and $\pi(x) \leq \pi(y)$, then
  $\pi(x) = \pi(x) \land \pi(y)$, so that the meet of $\pi(x),\pi(y)$
  exists in $Q$. Let $u \leq x,y$ be such that $\pi(u) = \pi(x) \land
  \pi(y) = \pi(x)$. Since $\pi$ is conservative, then $u = x$, so that
  $x \leq y$, exhibiting   $\pi$ as a Grothendieck fibration.

  Conversely, let us suppose that (i) and (ii) hold. Let $x,y \leq z$
  be such that $\pi(x) \land \pi(y)$ exists, and let $u,u' \leq x,y$
  be two preimages of this meet. Since $\pi$ is an embedding when
  restricted to the principal ideal of $z$ and $\pi(u) = \pi(u')$,
  then $u = u'$ as well. Let us show that $u = x \land y$: if $w \leq
  x, y$, then $\pi(w) \leq \pi(x),\pi(y)$ hence $\pi(w) \leq
  \pi(x)\land \pi(y) =\pi(u)$.  Since $w,u \leq z$, we deduce $w \leq
  u$, since $\pi$ is a Grothendieck fibration.
\end{proof}
An obvious remark, worth recalling at this point, is that if an order
preserving map $\pi : P \rTo Q$ creates pullbacks and $Q$ has
pullbacks, then $P$ has pullbacks as well which are preserved by
$\pi$. We state next the main result of this section.
\begin{theorem}
  \label{theo:charjsemidis}
  A lattice $L$ is join-semidistributive if and only if the projection
  $\pr[0] : \Cov(L) \rTo L$ creates pullbacks.
\end{theorem}
\begin{proof}
  If $\pr[0]$ creates pullbacks then, by the previous remark,
  $\Cov(L)$ has pullbacks. If $\mu \in \Cov(L)$ is maximal, then the
  ideal it generates is a finite meet-semilattice and hence it has a
  least element. Therefore $L$ is join-semidistributive.

  Conversely, let us suppose that $L$ is join-semidistributive.  Since
  $\pr[0]$ is conservative and, by Lemma \ref{lemma:joinsempushdown},
  a Grothendieck fibration, it is enough by Lemma
  \ref{lemma:picreatespbs} to show that if $\gamma, \delta\leq
  \epsilon$, then we can find $\beta \leq \gamma,\delta$ such that
  $\beta_{0} = \gamma_{0} \land \delta_{0}$.
  
  Observe first that $\gamma_{0} \land \delta_{0} < \gamma_{1}\land
  \delta_{1}$: otherwise, if $\gamma_{0} \land \delta_{0} =
  \gamma_{1}\land \delta_{1}$, then $\epsilon_{1} = \gamma_{1} \vee
  \epsilon_{0} = \delta_{1} \vee \epsilon_{0} = (\gamma_{1} \land
  \delta_{1}) \vee \epsilon_{0} = (\gamma_{0} \land \delta_{0}) \vee
  \epsilon_{0} = \epsilon_{0}$, a contradiction. 
  Let $\beta_{0} = \gamma_{0} \land \delta_{0}$ and choose $\beta_{1}$
  such that $\beta_{0} \lcover \beta_{1} \leq \gamma_{1}\land
  \delta_{1}$. We claim that $\beta_{1} \not\leq \gamma_{0}$:
  otherwise, if $\beta_{1} \leq \gamma_{0}$, then $\beta_{1} \leq
  \epsilon_{0} \land \delta_{1} = \delta_{0}$ and $\beta_{1} \leq
  \gamma_{0} \land \delta_{0} = \beta_{0}$ contradicting $\beta_{0}
  \lcover \beta_{1}$. 
  Therefore $\beta_{0} \leq \gamma_{0}$, $\beta_{1} \not\leq
  \gamma_{0}$, and $\beta_{1} \leq \gamma_{1}$, that is $\gamma \leq
  \beta$.  We derive $\beta \leq \delta$ similarly.
\end{proof}

The rest of this section is devoted to characterizing
join-semidistributive lattices among pushdown lattices.
\begin{lemma}
   \label{lemma:lowercover}
   If $L$ is a pushdown lattice and $\gamma \lcover \delta$ in the
   poset $\Cov(L)$, then $\delta \pushdown \gamma$.
   If $L$ is a join-semidistributive lattice and $\delta \pushdown
   \gamma$, then $\gamma \lcover \delta$ in $\Cov(L)$.
\end{lemma}
\begin{proof}
  Let us first assume  that $L$ is pushdown and 
  that $\gamma \lcover \delta$ in $\Cov(L)$. %
  Obviously, $\gamma < \delta$ and therefore 
  $\gamma_{1} < \delta_{1}$, by Lemma \ref{lemma:conservative}.  Let
  $u \in L$ be such that $\gamma_{1} \leq u \lcover \delta_{1}$ and
  observe that $u \neq \delta_{0}$, otherwise $\gamma_{1} \leq u =
  \delta_{0}$ contradicting $\gamma \leq \delta$. Let $\epsilon \in
  \Cov(L)$ be the unique cover such that $\delta \pushdown[u]
  \epsilon$. We have $\gamma \leq \delta$ and $\gamma_{0} \leq u$,
  hence $\gamma \leq \epsilon$ by Proposition
  \ref{prop:equivpushdown}.  
  Therefore $\gamma \leq \epsilon <
  \delta$, so that $\gamma \lcover \delta$ implies $\gamma =
  \epsilon$. 
  As $\gamma = \epsilon$, then we deduce that $\delta \pushdown[u]
  \gamma$.

  
  Next, we suppose that $L$ is a join-semidistributive lattice and
  that $\delta \pushdown[u] \gamma$. To prove that $\gamma \lcover
  \delta$, we consider $\epsilon \in \Cov(L)$ such that $\gamma \leq
  \epsilon \leq \delta$ and argue that either $\epsilon = \delta$, or
  $\epsilon \leq \gamma$.  If $\epsilon_{0} \not\leq u$, then
  $\epsilon_{1} \not\leq u$ and $\epsilon_{1} \vee u = \delta_{1}$. By
  join-semidistributivity, the relations $\epsilon_{1} \vee u =
  \delta_{1} = \epsilon_{1} \vee \delta_{0}$ imply $\delta_{1} =
  \epsilon_{1} \vee (u \land \delta_{0}) = \epsilon_{1} \vee
  \gamma_{0} = \epsilon_{1}$. By Lemma \ref{lemma:conservative},
  $\epsilon \leq \delta$ and $\epsilon_{1} = \delta_{1}$ imply
  $\epsilon = \delta$.
  Otherwise, if $\epsilon_{0} \leq u$, then Lemma
  \ref{lemma:joinsempushdown} and Proposition \ref{prop:equivpushdown}
  ensure that $\epsilon \leq \gamma$.
  %
\end{proof}

\begin{proposition}
  \label{prop:pushjsemid}
  For a pushdown lattice $L$ the following are equivalent:
  \begin{enumerate}
  \item If $\delta \pushdown \gamma$ then $\gamma \lcover \delta$,
  \item If $u\neq v$, $\delta \pushdown[u] \gamma$, and $\delta
    \pushdown[v] \epsilon$, then $\gamma_{0},\epsilon_{0}$ are not
    comparable.
  \item Every hat in $\Cov(L)$ has a pullback.
  \item $\Cov(L)$ has pullbacks.
  \item $L$ is join-semidistributive.
  \end{enumerate}
\end{proposition}
\begin{proof}
  \emph{(1) implies (2).} Let us suppose that $u\neq v$, $\delta \pushdown[u]
  \gamma$, and $\delta \pushdown[v] \epsilon$. 
  By the way of contradiction, let us suppose that $\gamma_{0} \leq
  \epsilon_{0}$. Then $\gamma, \epsilon < \delta$ and $\gamma_{0} \leq
  \epsilon_{0}$ imply $\gamma \leq \epsilon$. Taking into account
  that, by (1), $\gamma \lcover \delta$, then $\gamma \leq \epsilon <
  \delta$ implies $\gamma = \epsilon$. Consequently, $\delta
  \pushdown[u] \gamma$ and $\delta \pushdown[v] \gamma$ imply
  $\gamma_{1} \leq u \land v$.
  
  Next, let $\delta' = (u,\delta_{1})$ and let $\psi,\xi \in \Cov(L)$
  be the covers determined by the relations $\delta'
  \spushdown[3pt]{\delta_{0}} \psi$ and $\delta' \pushdown[v] \xi$.
  Observe that $\psi_{0} < \xi_{0}$, since $\psi_{0} = u \land
  \delta_{0} = \gamma_{0} < \gamma_{1} \leq u \land v = \xi_{0}$.
  Thus, considering that $\psi_{0} < \xi_{0}$ and $\psi,\xi \leq
  \delta'$, we can use the pushdown property 
  to deduce that $\psi < \xi$.  However, this is a contradiction: we
  have $\psi < \xi < \delta'$ and, by (1), $\psi \lcover \delta'$.


  \breath
  \noindent
  \emph{(2) implies (3).}  Let us suppose that $\gamma ,\epsilon\lcover
  \delta$ with $\gamma \neq \epsilon$. By Lemma \ref{lemma:lowercover}
  we can write $\delta \pushdown[u] \gamma$ and $\delta \pushdown[v]
  \epsilon$ for some $u,v \in L$ such that, by Lemma
  \ref{lemma:uniquepusher}, $u \neq v$.  Hence, $u,v,\delta_{0}$ are
  pairwise distinct lower covers of $\delta_{1}$.
  
  We claim that, assuming (2), $\gamma_{0} \land \epsilon_{0} < u\land
  v$. Clearly, we have $\gamma_{0} \land \epsilon_{0} \leq u\land v$,
  so that, by the way of contradiction, let us assume that $\gamma_{0}
  \land \epsilon_{0} = u\land v$; it follows that $u\land v \leq
  \gamma_{0}$. Let $\delta' = (u,\delta_{1})$ and let $\psi,\xi \in
  \Cov(L)$ be the covers determined by $\delta'
  \spushdown[3pt]{\delta_{0}} \psi$ and $\delta' \pushdown[v] \xi$.
  It follows that $\xi_{0} \leq \psi_{0}$, since $\xi_{0} =u \land v
  \leq \gamma_{0} = u \land \delta_{0} = \psi_{0}$, thus contradicting
  property (2).

  Next, we construct a lower bound $\beta$ of $\gamma,\epsilon$. 
  Given
  that $\gamma_{0} \land \epsilon_{0} < u\land v$, we let $\beta_{0}=
  \gamma_{0} \land \epsilon_{0}$ and pick $\beta_{1}$ such that
  $\beta_{0} \lcover \beta_{1} < u \land v$. 
  Let us argue that $\beta
  = (\beta_{0},\beta_{1})$ is below $\gamma$.
  Clearly, $\beta_{0} = \gamma_{0} \land \epsilon_{0} \leq
  \gamma_{0}$. Consequently $\beta_{0} \leq \delta_{0}$ and, by
  construction, $\beta_{1} \leq u \land v \leq \delta_{1}$. If
  $\beta_{1} \leq \delta_{0}$ then $\beta_{1} \leq \delta_{0} \land u
  \land v = (\delta_{0} \land u) \land (\delta_{0} \land v) =
  \gamma_{0} \land \epsilon_{0} = \beta_{0}$, a contradiction; whence
  $\beta_{1} \not\leq \delta_{0}$. We have argued that $\beta \leq
  \delta$ and $\beta_{0} \leq \gamma_{0}$ which, together with $\gamma
  \leq \delta$, imply $\beta \leq \gamma$. 
  In a similar way we show that $\beta \leq \epsilon$, so that $\beta$
  is a lower bound of $\gamma,\epsilon$.

  Finally, to see that $\beta$ is the greatest lower bound of $\gamma$
  and $\epsilon$, we argue as in Lemma \ref{lemma:picreatespbs}. If
  $\alpha \leq \gamma, \epsilon$ then $\alpha_{0} \leq \gamma_{0}
  \land \epsilon_{0} = \beta_{0}$ and, considering that $\alpha,\beta
  \leq \delta$ and $L$ is pushdown, we deduce $\alpha \leq \beta$.
  

  \breath
  \noindent
  \emph{(3) implies (4).} By Lemma   \ref{lemma:pullbacksofhats}.

  \breath
  \noindent
  \emph{(4) implies (5).} %
  Since $\Cov(L)$ has pullbacks, if $\mu \in \Cov(L)$ is a maximal
  cover, then the ideal it generates is a finite meet-semilattice and
  has a least element.  By Lemma \ref{lemma:usual}, it follows that
  $L$ is join-semidistributive.
  
  

  \breath
  \noindent
  \emph{(5) implies (1).} By Lemma \ref{lemma:lowercover}.
\end{proof}


%% file: bounded.tex
\section{Lower Bounded Lattices}
\label{sec:lbounded}
\label{sec:bounded}

Recall that a lattice $L$ is said to be \emph{lower bounded} if there
exists a finite set $X$ and lattice epimorphism from the freely
generated lattice ${\cal F}(X)$\footnote{If $X$ contains at least
  three elements, then ${\cal F}(X)$ is an infinite lattice.}  onto
$L$, say $f : {\cal F}(X) \rOnto L$, such that, for each $y \in L$,
the set $\set{ x \mid y \leq f(x) }$ is either empty or has a least
element, see \cite[\S 5]{mckenzie} or \cite[\S
2.13]{freese}. \emph{Upper boundedness} is the dual notion of lower
boundedness, and a lattice is said to be \emph{bounded} if it is both
lower and upper bounded.

There are already many characterizations of  lower bounded
lattices \cite{day,jonssonnation2,urquhart,semenova} and this concept
has also found applications within unexpected branches of lattice
theory
\cite{gratzerwerhung}.  In this section we develop further the tools
used in \cite{BCM} to prove that lattices in ${\cal H}{\cal H}$ -- a
class which generalizes the class of finite Coxeter lattices -- are
bounded. In this way we shall obtain a new characterization of lower
bounded lattices.  Our starting point will be the following classical
result \cite[\S 2.39]{freese}:
\begin{theorem}[Johnsonn, Nation] 
  \label{theo:Dnocycles}
  A lattice is lower bounded if and only if the \emph{join-dependency
    relation} between join-irreducible elements contains no cycle.
\end{theorem}
In order to use of this result, we must of course give the
definition of the join-dependency relation $D$.
Here we give its characterization in terms of
the arrows relations \cite[\S 11.10]{freese}.
The reader can find
in \cite[\S 2.3]{freese} its standard definition. 
\begin{definition}
  For a lattice $L$ and $j,k \in J(L)$, we let $j D k$ if $j \neq k$ and
  there exists $m \in M(L)$ such that $j \jmup m$ and $k \jmdown m$.
\end{definition}
Let us introduce some more relations between join-irreducible elements
of a lattice $L$:
\begin{definition}
  For all $j,k \in J(L)$, we let:
  \begin{itemize}
  \item $j A k$ if $j \neq k$ and, for some $m \in M(L)$, $j \jmup m$
    and $k \jmpersp m$,
  \item $j B k$ if $j \neq k$ and, for some $m \in M(L)$,
    $j \jmpersp m$ and $k \jmdown m$,  
  \item $j C k$ if either $j A k$ or $j B k$.
  \end{itemize}
\end{definition}
Let us remark that these relations are already known for
semidistributive lattices \cite[\S 2.5]{freese}. However, the
definition presented here makes sense in any lattice. The next Lemma
shows that, in the expression of Theorem \ref{theo:Dnocycles}, the $D$
relation may be replaced by the $C$ relation defined above.
\begin{lemma}
  \label{lemma:replaceDwC}
  The $D$ relation has a cycle if and only if the $C$ relation has a
  cycle. Consequently a  lattice is lower bounded if and only if
  the $C$ relation contains no cycle.
\end{lemma}
\begin{proof}
  Let us suppose that $j D k$ and let $m \in M(L)$ be such that $j
  \jmup m$ and $k \jmdown m$. Choose $l \in J(L)$ such that $l
  \jmpersp m$.  If $l \in \set{j,k}$, then $j A k$ or $j B k$.  If
  $l\neq j,k$, then $j A l$ and $l B k$. Therefore one step of the
  relation $D$ may be replaced by at most two steps of the relation
  $C$.  Conversely, since $C \subseteq D$, every $C$-cycle gives rise
  to a $D$-cycle.
\end{proof}


We shall introduce a number of relations and, to this goal, we first
need the dual of the relation $\pushdown$: for $w \in L$ and
$\gamma,\delta \in \Cov(L)$, let us write $\gamma \pushup[w] \delta$
if $w \neq \gamma_{1}$, $\delta_{1} = w \vee \gamma_{1}$, and
$\gamma_{0} \lcover w \leq \delta_{0}$. Also let us write $\gamma
\pushup \delta$ if $\gamma \pushup[w] \delta$ for some $w \in L$.
The following remark  is due now:
\begin{lemma}
  If a lattice $L$ is join-semidistributive, then $\delta \pushdown
  \gamma$ implies $\gamma \pushup \delta$.
  Consequently, 
  if $L$ is a semidistributive lattice, then $\delta \pushdown \gamma$
  holds if and only if $\gamma \pushup \delta$ holds.
  \label{lemma:pushdownpushup}
\end{lemma}
\begin{proof}
  Let us suppose that $\delta \pushdown[u] \gamma$ so that, as
  usual, $\gamma_{0} = u \land \delta_{0} < \delta_{0}$. Thus, we can
  choose $w \in L$ such that $\gamma_{0} \lcover w \leq
  \delta_{0}$. Forcedly, $w \not\leq u$, since otherwise $\gamma_{0} <
  w \leq u \land \delta_{0} = \gamma_{0}$.  It follows that $w \neq
  \gamma_{1}$ and that $w \vee u = \delta_{1}$.  We claim that $w \vee
  \gamma_{1} = \delta_{1}$.  If this is not the case, then $w \vee
  \gamma_{1} < \delta_{1}$ and the implication \eqref{eq:sdvee} fails:
  let $x = w \vee \gamma_{1}$, $y = u$, and $z = \delta_{0}$, then $x
  \vee y = w \vee \gamma_{1} \vee u = w \vee u = \delta_{1} =
  \gamma_{1} \vee \delta_{0} = w \vee \gamma_{1} \vee \delta_{0} = x
  \vee z$, but $x \vee (y \land z) = x \vee \gamma_{0} = x <
  \delta_{1} = x \vee y$.
\end{proof}


\begin{definition}
  For $\gamma,\delta \in \Cov(L)$, we let $\gamma \cAm \delta$ if there
  exists $\epsilon \in \Cov(L)$ and $u \in L$ such that $\epsilon
  \pushdown[u] \delta$ and $\delta_{1} \leq \gamma_{0} \lcover
  \gamma_{1} \leq u$. The dual relation is defined as follows: $\gamma
  \cBm \delta$ if there exists $\epsilon \in \Cov(L)$ and
  $u \in L$ such that $\epsilon \pushup[u] \delta$ and $u \leq
  \gamma_{0} \lcover \gamma_{1} \leq \delta_{1}$.
\end{definition}
Intuitively, the $\cAm,\cBm$ relations express the dependency relation
of covers in terms of the relations $\pushdown$ and $\pushup$, when
the configurations they give rise to form pentagons.  The following
picture illustrates this point with the $\cAm$ relation:
$$
\xygraph{
  []!~:{@{.}}
  !d{0.5}{\epsilon}
  "\epsilon_{1}"-[d(0.7)l]*+{u}="u"
  [d(1.5)]!c{\gamma}"u":"\gamma_{1}"
   "\gamma_{0}"[d(1.5)r(0.5)]!d{0.5}{\delta}
  "\gamma_{0}":"\delta_{1}"
  "\epsilon_{0}":"\delta_{0}"
  "\gamma_{m}"[r(0.15)]="\gamma_{m}"
  "\delta_{m}"[r(0.1)u(0.1)]="\delta_{m}"
  "\gamma_{m}":@{{ }{-}{>}}@/^1.0em/"\delta_{m}"^{\cAm}
}
$$
It would be tempting to define \emph{facets} at this
point. Intuitively, these are the configurations arising either from a
relation of the form $\delta \pushdown[u] \gamma$, or from a relation
of the form $\gamma \pushup[w] \delta$. We believe, however, that the
notion of facet is better suited to the semidistributive case, where
the two relations $\pushdown$ and $\pushup$ are one opposite of the
other, as remarked in Lemma \ref{lemma:pushdownpushup}. Since we want
to carry out our ideas in the weaker join-semidistributive setting, we
delay the definition of facets until the end of this section when we
shall come back to full semidistributivity.


\medskip

The next two Lemmas exemplify the connections between the relations $A$
and $\cAm$. The Lemmas apply to arbitrary  lattices.
\begin{lemma}
  \label{lemma:AAtoA}
  Let $j, k \in J(L)$ be such that, for some $\gamma, \delta \in
  \Cov(L)$, $\gamma \cAm \delta$, $(j_{\ast},j) \leq \gamma$ and
  $(k_{\ast},k) \leq \delta$.  Then $j A k$.
\end{lemma}
\begin{proof}
  Let $\gamma,\delta$ be such that $\gamma \cAm \delta$ so that
  $\epsilon \pushdown[u] \delta$, for some $u \in L$ and $\epsilon \in
  \Cov(L)$, and moreover $\delta_{1} \leq \gamma_{0} \lcover
  \gamma_{1} \leq u$.  Observe first that $j \neq k$, since $k \leq
  \delta_{1} \leq \gamma_{0}$ and $j \not\leq \gamma_{0}$.
  Let $\mu$ be maximal in the set $\set{ \theta \in \Cov(L) \mid
    \epsilon \leq \theta}$ so that $\mu$ is of the form $\mu =
  (m,m^{\ast})$ for some $m \in M(L)$.
  We claim that $j \jmup m$.  We have on one side $j \leq \gamma_{1}
  \leq u \leq \epsilon_{1} \leq m^{\ast}$. Let us suppose that $j \leq
  m$: then $j \leq u \land \epsilon_{1} \land m = u \land \epsilon_{0}
  = \delta_{0} \leq \gamma_{0}$. This is a contradiction, since
  $(j_{\ast},j) \leq \gamma$ implies $j \not \leq \gamma_{0}$;
  hence $j \not\leq m$. 
  We have shown that $j \jmup m$; since
  $(k_{\ast},k) \leq \delta \leq \epsilon \leq (m,m^{\ast})$, we have
  $k \jmpersp m$, and therefore $j A k$.
\end{proof}

In the next Lemma~$\,\pushdown[\ast]$ denotes the reflexive transitive
closure of the pushdown relation $\pushdown$. Since $\delta
\pushdown[\ast] \gamma$ implies $\gamma \leq \delta$, the relation
$\pushdown[\ast]$ is antisymmetric, and hence it is an ordering on
$\Cov(L)$. We recall that this ordering coincides with the opposite of
the order $\leq$ on $\Cov(L)$ if $L$ is a pushdown lattice, see
Proposition \ref{prop:transcolos}; on the other hand, the scope of the
next Lemma is not restricted to pushdown lattices.
\begin{lemma}
  \label{lemma:pushdA}
  Let $j \in J(L)$ and $\delta \in \Cov(L)$ be such that $\delta_{1} =
  j \vee \delta_{0}$.  Then there exists $n \geq 0$ and a sequence
  $(\delta^{i},\gamma^{i})$, $i = 0,\ldots ,n$, such that 
  \begin{enumerate}
  \item $\gamma^{0} = (j_{\ast},j)$ and $\delta^{n} = \delta$,
  \item $\delta^{i}_{1} = j \vee \delta^{i}_{0}$ and $\delta^{i}
    \pushdown[\ast] \gamma^{i}$,
    for $i = 0,\ldots,n$,
  \item $\delta^{i-1} \cAm \gamma^{i}$, for $i = 1,\ldots ,n$.
  \end{enumerate}
\end{lemma}
\begin{proof}
  The proof is by induction on the height of the interval
  $[j,\delta_{1}]$. 

  Consider the set
  \begin{align*}
    PD(\delta,j)& = \set{\theta \in \Cov(L) \mid \delta
      \pushdown[\ast] \theta \tand j \vee \theta_{0} = \theta_{1}}\,,
  \end{align*}
  observe that $\delta \in PD(\delta,j)$, and hence pick $\mu \in
  PD(\delta,j)$ which is maximal w.r.t. the relation
  $\pushdown[\ast]$.

  If $j = \mu_{1}$, then $\mu_{1}$ is join-irreducible so that
  $\mu_{0} = j_{\ast}$ and the statement holds with $n = 0$,
  $\delta^{0} = \delta$, and $\gamma^{0} = \mu = (j_{\ast},j)$.

  If $j \neq \mu_{1}$, then $j < \mu_{1}$ and we can find $u
  \in L$ such that $j \leq u \lcover \mu_{1}$.  As usual, observe
  that $u \neq \mu_{0}$, since $j \not\leq \mu_{0}$.
  We let therefore $\gamma_{0} = u \land \mu_{0}$ and $\delta'_{1} = j
  \vee \gamma_{0}$.  We observe next two facts about the pair
  $(\gamma_{0} ,\delta'_{1})$: (a) we have $\gamma_{0} < \delta'_{1}$,
  since if $\gamma_{0} = \delta'_{1} = j \vee \gamma_{0}$, then $j
  \leq \gamma_{0} \leq \mu_{0}$, contradicting $j \vee \mu_{0} =
  \mu_{1}$; (b) the pair $(\gamma_{0},\delta'_{1})$ is not a cover,
  since otherwise, defining $\delta' = (\gamma_{0},\delta'_{1})$, then
  $\mu \pushdown[u] \delta'$ and $\delta'_{1} = j \vee \gamma_{0}$
  imply that $\mu$ is not maximal in the set $PD(\delta,j)$.  As a
  consequence of our observations, we can find $\gamma_{1},
  \delta'_{0} \in L$ such that $\gamma_{0} \lcover \gamma_{1} \leq
  \delta'_{0} \lcover \delta'_{1} \leq u$. If we let $\gamma =
  (\gamma_{0},\gamma_{1})$ and $\delta' = (\delta'_{0},\delta'_{1})$,
  then $\mu \pushdown[u] \gamma$ and $\delta'\cAm \gamma$;
  consequently, $\delta \pushdown[\ast] \gamma$ as well.


  Next we remark that $[j,\delta'_{1}] \subset [j,\delta_{1}]$, 
  since  $\delta'_{1} \leq u < \delta_{1}$.
  We can use the inductive hypothesis to find $m \geq 0$ and a
  sequence $(\delta^{i},\gamma^{i})$, $i = 0,\ldots m$, satisfying
  \enumref{2} and \enumref{3} and such that $\gamma^{0} = (j_{\ast},j)$ and
  $\delta^{m} = \delta'$.
  %
  %
  Therefore we let $n = m+1$ and, to obtain a sequence satisfying
  \enumref{1}, \enumref{2}, and \enumref{3}, we append the pair
  $(\delta,\gamma)$ to the sequence $(\gamma^{i},\delta^{i})$, $i =
  0,\ldots ,m$.
\end{proof}

Following \cite{BCM}, we are ready to introduce strict facet
labellings. 
\begin{definition}
  \label{def:labelling}
  \label{def:facetlabelling}
  A \emph{strict lower facet labelling} of a lattice $L$ is a function
  $f : \Cov(L) \rTo \N$ such that $f(\gamma) = f(\delta)$ if $\delta
  \pushdown\gamma$ and $f(\gamma) < f(\delta)$ if $\delta \cAm
  \gamma$.
  A \emph{strict upper facet labelling} is defined dually: it is a
  function $f : \Cov(L) \rTo \N$ such that $f(\gamma) = f(\delta)$ if
  $\delta \pushup \gamma$ and $f(\gamma) < f(\delta)$ if $\delta \cBm
  \gamma$.

  
  A \emph{strict facet labelling} of a lattice $L$ is a function $f :
  \Cov(L) \rTo \N$ which is both a strict lower facet labelling and a
  strict upper facet labelling.
\end{definition}
\begin{lemma}
  \label{lemma:equality}
  If $L$ is a pushdown lattice, $f$ is a strict lower facet labelling,
  and $\gamma \leq \delta$, then $f(\gamma) = f(\delta)$.
\end{lemma}
\begin{proof}
  If $\gamma \leq \delta$ then by Lemma \ref{lemma:lowercover} we can
  find a path of the relation $\pushdown$ from $\delta$ to $\gamma$.
  The statement follows since $f$ is constant on the relation
  $\pushdown$.
\end{proof}

\begin{lemma}
  \label{lemma:downstr}
  Let $L$ be a lattice with a strict lower facet labelling $f$.  Let
  $j \in J(L)$ and $\delta \in \Cov(L)$ be such that $\delta_{1} = j
  \vee \delta_{0}$. Then $f(\delta) \leq f(j_{\ast},j)$ and $f(\delta)
  = f(j_{\ast},j)$ implies $(j_{\ast},j) \leq \delta$.
\end{lemma}
\begin{proof}
  By Lemma \ref{lemma:pushdA}, let $n \geq 0$ and
  $(\delta^{i},\gamma^{i})$, $i = 0,\ldots ,n$, be such that
  $\delta^{i} \pushdown[\ast] \gamma^{i}$ for $i = 0,\ldots ,n$,
  $\delta^{i -1}\cAm \gamma^{i}$ for $i = 1,\ldots ,n$, $\delta^{n} =
  \delta$ and $\gamma^{0} = (j,j_{\ast})$.

  Since $f$ is constant on $\pushdown$, it is constant on its
  reflexive transitive closure $\pushdown[\ast]$. Hence, we have
  $f(\delta^{i}) = f(\gamma^{i})$ for $i = 0,\ldots ,n$, and
  $f(\gamma^{i}) < f(\delta^{i-1})$ for $i = 1,\ldots ,n$.  We deduce
  that $f(\delta) = f(\delta^{n})\leq f(\gamma^{0}) = f(j_{\ast},j)$
  and, $f(\delta) < f(j_{\ast},j)$ if $n \geq 1$.  Therefore, if
  $f(\delta) = f(j_{\ast},j)$, then $n =0$, $\delta \pushdown[\ast]
  (j_{\ast},j)$ and $(j_{\ast},j) \leq \delta$.
\end{proof}

By duality, we obtain the following Corollary.
\begin{corollary}
  \label{lemma:upstr}
  Let $L$ be a 
  lattice with an upper strict facet labelling
  $f$.  Let $m \in M(L)$ and $\delta \in \Cov(L)$ such that $m \land
  \delta_{1} = \delta_{0}$.  Then $f(\delta) \leq f(m,m^{\ast})$ and
  $f(\delta) = f(m,m^{\ast})$ implies $\delta \leq (m,m^{\ast})$.
\end{corollary}

We are ready to achieve the first goal of this section, Proposition
\ref{prop:facettobounded}.
\begin{lemma}
  \label{lemma:A}
  If $L$ is a join-semidistributive lattice with a strict lower facet
  labelling $f$, then $j A k$ implies $f(k_{\ast},k) < f(j_{\ast},j)$.
\end{lemma}
\begin{proof}
  If $j A k$ then $j \neq k$, $j \jmup m$, and $(k_{\ast},k) \leq (m,m^{\ast})$ for
  some $m \in M(L)$. It follows from Lemma \ref{lemma:equality} that
  $f(k_{\ast},k) = f(m,m^{\ast})$.
  Since moreover we have $j \vee m = m^{\ast}$, Lemma
  \ref{lemma:downstr} implies that $f(m,m^{\ast}) \leq f(j_{\ast},j)$
  and moreover that this is an inequality: otherwise, if
  $f(m,m^{\ast}) = f(j_{\ast},j)$, then $(j_{\ast},j) \leq
  (m,m^{\ast})$; consequently, $j = k$ since in a
  join-semidistributive lattice there exists at most one $j \in
  \Cov(L)$ such that $j \jmpersp m$.
\end{proof}

The next Lemma is not a mere consequence of duality.
\begin{lemma}
  \label{lemma:B}
  If $L$ is a join-semidistributive lattice with a strict upper
  facet labelling $f$, then $j B k$ implies $f(k_{\ast},k) <
  f(j_{\ast},j)$.
\end{lemma}
\begin{proof}
  Let $j,k \in J(L)$ be such that $jB k$, that is, $j \neq k$ and, for
  some $m \in M(L)$, $j \jmpersp m$ and $k \jmdown m$.  Since by Lemma
  \ref{lemma:joinsempushdown} $L$ is a pushup lattice, we can make use
  of Lemma \ref{lemma:equality} to deduce that $f(j_{\ast},j) =
  f(m,m^{\ast})$. On the other hand, Lemma \ref{lemma:upstr} ensures
  that $f(k_{\ast},k) \leq f(m,m^{\ast})$ and $(k_{\ast},k) \leq
  (m,m^{\ast})$ if $f(k_{\ast},k) = f(m,m^{\ast})$. We cannot have
  $f(k_{\ast},k) = f(m,m^{\ast})$, since then $(k_{\ast},k) \leq
  (m,m^{\ast})$, $k \jmpersp m$, and $j = k$ by
  join-semidistributivity. Hence $f(k_{\ast},k) < f(m,m^{\ast}) =
  f(j_{\ast},j)$.
\end{proof}

  

\begin{proposition}
  \label{prop:facettobounded}
  A join-semidistributive lattice with a strict facet labelling is
  lower bounded.
\end{proposition}
\begin{proof}
  According to Lemma \ref{lemma:replaceDwC}, it is enough to show that
  the relation $C$ has no cycle, for which we shall argue that $j C k$
  implies $f(k_{\ast},k) < f(j_{\ast},j)$.  If $j A k$ then we use
  Lemma \ref{lemma:A}. If $j B k$ then we can use Lemma \ref{lemma:B}.
\end{proof}

Our next goal is to prove the converse of Proposition
\ref{prop:facettobounded}. 


\begin{lemma}
  \label{lemma:BBtoB}
  Let $L$ be a join-semidistributive lattice.  If $j, k \in J(L)$ are
  such that for some $\gamma, \delta \in \Cov(L)$, $\gamma \cBm
  \delta$, $(j_{\ast},j) \leq \gamma$ and $(k_{\ast},k) \leq \delta$,
  then $j B k$.
\end{lemma}
\begin{proof}
  Let $\gamma,\delta$ be such that $\gamma \cBm \delta$: for some $u
  \in L$ and $\epsilon \in \Cov(L)$, $\epsilon \pushup[u] \delta$ and
  $\epsilon_{0} \lcover u \leq \gamma_{0}\lcover \gamma_{1} \leq
  \delta_{0}$. Let $j,k$ be as in the statement, and observe that $j
  \neq k$, since $j \leq \gamma_{1} \leq \delta_{0}$ and $k \not\leq
  \delta_{0}$.
  Let $\mu = (m,m^{\ast})$ be maximal above $\gamma$, so that $j
  \jmpersp m$, and let $l \in J(L)$ such $(l_{\ast},l) \leq \epsilon$.
  We have $m \geq \gamma_{0} \geq u \geq \epsilon_{0} \geq l_{\ast}$.
  If $m \geq l$, then $m \geq u \vee \epsilon_{0} \vee l = u \vee
  \epsilon_{1} = \delta_{1}$; this is in turn implies that $m \geq
  \delta_{1} \geq \gamma_{1}$, a contradiction. Hence $m \not\geq l$
  and $l \jmdown m$. We have therefore $j B l$.
  It is now easy to see that $l = k$, since $(l_{\ast},l) \leq \delta$
  and $(k_{\ast},k) \leq \delta$ imply $l = k$, by join-semidistributivity.
\end{proof}

\begin{proposition}
  A lower bounded lattice has a strict facet labelling.
\end{proposition}
\begin{proof}
  If $L$ is bounded then it is join-semidistributive \cite[\S
  2.20]{freese} and the
  join-dependency relation $D$ is acyclic \cite[\S 2.39]{freese}.
  Let us denote by $\Dclosure$ its reflexive and transitive closure,
  and let $g : J(L)\rTo \N$ be an antilinear extension of the poset
  $\langle J(L),\Dclosure\rangle$. That is, $g$ is such that $g(y) <
  g(x)$ whenever $x D y$.

  We define $f : \Cov(L)\rTo \N$ as follows:
  $f(\delta) = g(j(\delta))$ where $j(\delta) \in J(L)$ is the unique
  join-irreducible element such that $(j(\delta)_{\ast},j(\delta))
  \leq \delta$.  Then $f$ is a strict facet labelling. It is a lower
  facet labelling: if $\delta \pushdown \gamma$ then $j(\delta) =
  j(\gamma)$ and $f(\delta) = g(j(\delta)) = g(j(\gamma)) =
  f(\gamma)$; if $\delta \cAm \gamma$, then $j(\delta) A j(\gamma)$,
  by Lemma \ref{lemma:AAtoA}, and $j(\delta) D j(\gamma)$; it follows
  that $f(\gamma) = g(j(\gamma)) < g(j(\delta)) = f(\delta)$.
  Similarly, it is an upper strict facet labelling, since $\delta
  \pushup[u] \gamma$ implies $j(\gamma) = j(\delta)$ and $\delta \cBm
  \gamma$ implies $j(\delta) B j(\gamma)$, by Lemma \ref{lemma:BBtoB},
  and therefore $f(\gamma) < f(\delta)$.
\end{proof}

Recalling that lower bounded lattices are join-semidistributive
\cite[\S 2.20]{freese}, we end this section collecting the
observations presented so far into a main result:
\begin{theorem}
  \label{theo:mainresbounded}
  A  lattice is lower bounded if and only if it is
  join-semidistributive and has a strict facet labelling.
\end{theorem}
Finally, we observe that from Theorem \ref{theo:mainresbounded} it is
quite immediate to derive a standard result by Day, see \cite[\S
2.64]{freese}, stating that a semidistributive lower bounded lattice
is bounded.  Recall now that in a semidistributive lattice $\delta
\pushdown \gamma$ holds if and only if $\gamma \pushup \delta$ holds.
Theorem \eqref{theo:mainresbounded} leads to a characterization of
bounded lattices, that we shall rephrase in a language closer to
\cite{BCM}. Call a \emph{facet} of $L$ a quadruple of distinct covers
$(\delta,\delta',\gamma,\gamma') \in \Cov^{4}(L)$ such that
\begin{align*}
\delta_{1} = \delta'_{1} & \tand \gamma_{0} =
  \gamma'_{0}\,,
  \text{ and moreover }
  \delta \spushdown[4pt]{\delta'_{0}} \gamma  \tand
  \delta' \spushdown[3pt]{\delta_{0}} \gamma'\,.
\end{align*}
See Figure 1 in the Introduction.  Say that $\epsilon$ is
\emph{interior} to the facet $(\delta,\delta',\gamma,\gamma')$ if
$\gamma_{1} \leq \epsilon_{0} \lcover \epsilon_{1} \leq \delta'_{0}$
or $\gamma'_{1} \leq \epsilon_{0} \lcover \epsilon_{1} \leq
\delta_{0}$.
\begin{theorem}
  \label{theo:bounded}
  A lattice is bounded if and only if it is semidistributive and there
  exists a function $f: \Cov(L) \rTo \N$ such that for each facet 
  $(\delta,\delta',\gamma,\gamma')$ $f(\delta) =
  f(\gamma)$ and $f(\delta') = f(\gamma')$, and moreover
  $f(\delta'),f(\gamma) < f(\epsilon)$ whenever $\epsilon$ is
  interior to such a facet.  
\end{theorem}


%% file: derived.tex
\section{Derived Semidistributive Lattices}
\label{sec:derived}

The main result of this section is Theorem \ref{theo:derived} stating
that a poset of the form $\Cov(L,\gamma)$, $\gamma \in \Cov(L)$, is a
semidistributive lattice whenever $L$ is semidistributive.  Observe
that if $\gamma = (j_{\ast},j)$ with $j \in J(L)$ and $L$ is
join-semidistributive, then $\Cov(L,\gamma)$ is the set $\set{\delta
  \in \Cov(L)\mid\gamma \leq \delta}$. With this in mind we observe:
\begin{proposition}
  If $L$ is a semidistributive lattice, then $\Cov(L,\gamma)$ is a
  lattice, for each $\gamma \in \Cov(L)$.
\end{proposition}
\begin{proof}
  By Proposition \ref{prop:pushjsemid} and its dual $\Cov(L)$ has
  pullbacks and pushouts, hence $\Cov(L,\gamma)$ is a lattice by
  Corollary \ref{cor:pullpush}.
\end{proof}
We shall call $\Cov(L,\gamma)$ the semidistributive lattice
\emph{derived} from $L$ and $\gamma$. We study next additional
properties of the lattices of the form $\Cov(L,\gamma)$.
The next Lemma will prove useful in establishing these properties.
\begin{lemma}
  \label{lemma:useful}
  Let $L$ be a join-semidistributive lattice, and let
  $\gamma,\delta,\epsilon \in \Cov(L)$ be such that $\gamma \leq \delta$,
  $\epsilon \lcover \delta$, and $\gamma \not\leq \epsilon$.  Then
  $\gamma_{0} \vee \epsilon_{1} = \delta_{1}$.
\end{lemma}
\begin{proof}
  From $\gamma \leq \delta$ and $\epsilon \lcover \delta$, it follows
  that $\gamma_{1}\vee \epsilon_{1} \leq \delta_{1}$. We claim that
  $\gamma_{1}\vee \epsilon_{1} = \delta_{1}$. Observe that if the
  claim holds, then we can easily derive the conclusion of the Lemma,
  $\gamma_{0} \vee \epsilon_{1} = \delta_{1}$, since by
  join-semidistributivity $ \gamma_{1} \vee \epsilon_{1} = \delta_{1}
  = \delta_{0} \vee \epsilon_{1}$ implies $ \gamma_{0} \vee
  \epsilon_{1} = (\gamma_{1} \land \delta_{0}) \vee \epsilon_{1} =
  \gamma_{1} \vee \epsilon_{1} = \delta_{1}$.

  Let us suppose next that the claim does not hold, i.e. that
  $\gamma_{1}\vee \epsilon_{1} < \delta_{1}$.  Choose $u \in L$ such
  that $\gamma_{1} \vee \epsilon_{1} \leq u \lcover \delta_{1}$ and
  observe that $u \neq \delta_{0}$, since otherwise $\gamma,\epsilon
  \not\leq \delta$.  Let $\epsilon' \in \Cov(L)$ be such that $\delta
  \pushdown[u] \epsilon'$.  Since $\gamma_{0}\leq \gamma_{1} \leq u$
  and $\epsilon_{0} \leq \epsilon_{1} \leq u$, Proposition
  \ref{prop:equivpushdown} implies that $\gamma \leq \epsilon'$ and
  $\epsilon \leq \epsilon'$. In turn, $\epsilon \leq \epsilon' <
  \delta$ and $\epsilon \lcover \delta$ imply $\epsilon =
  \epsilon'$. We have therefore $\gamma \leq \epsilon' = \epsilon$,
  contradicting the hypothesis of the Lemma, $\gamma \not\leq
  \epsilon$.
%
\end{proof}

In the following, we shall use capital Greek letters to range on
elements of $\Cov(\Cov(L,\gamma))$.  Observe that the next
Propositions make sense: if $L$ is join-semidistributive then by
Proposition \ref{prop:pushjsemid} $\Cov(L,\gamma)$ is a finite poset
with pullbacks. This is enough to ensure that the relation $\leq$ on
$\Cov(\Cov(L,\gamma))$ is a transitive relation. Since this relation
is clearly reflexive and antisymmetric, it is a partial ordering.
\begin{proposition}
  If $L$ if a join-semidistributive lattice then the projection
  \begin{align*}
    \pr[0] : & \,\Cov(\Cov(L,\gamma)) \rTo \Cov(L,\gamma) 
  \end{align*}
  is a Grothendieck fibration.
\end{proposition}
\begin{proof}
  Let $\Gamma,\Psi,\Delta \in \Cov(\Cov(L,\gamma))$ be such that
  $\Gamma,\Psi \leq \Delta$ and $\Gamma_{0} \leq \Psi_{0}$. We shall
  prove that $\Gamma_{1} \leq \Psi_{1}$, which in turn  implies
  that $\Gamma \leq \Psi$.

  Observe that, since $\Gamma_{1},\Psi_{1} \leq \Delta_{1}$, the
  pullback $\Gamma_{1} \land \Psi_{1}$ exists in $\Cov(L)$. Moreover,
  since $\Gamma_{0} \leq \Psi_{0} \leq \Psi_{1}$, we have $\Gamma_{0}
  \leq \Gamma_{1} \land \Psi_{1} \leq \Gamma_{1}$ and therefore either
  $\Gamma_{1} \land \Psi_{1} = \Gamma_{1}$, or $\Gamma_{1} \land
  \Psi_{1} = \Gamma_{0}$, since $\Gamma_{0} \lcover \Gamma_{1}$.  We
  shall show that $\Gamma_{1} \land \Psi_{1} = \Gamma_{0}$ cannot
  occur, hence we have $\Gamma_{1} \land \Psi_{1} = \Gamma_{1}$, that
  is, $\Gamma_{1} \leq \Psi_{1}$.

  By the way of contradiction, assume that $\Gamma_{1} \land \Psi_{1}
  = \Gamma_{0}$.  Recall from Theorem \ref{theo:charjsemidis} that
  $\Gamma_{1,0} \land \Psi_{1,0} = (\Gamma_{1} \land \Psi_{1})_{0}$,
  and hence we see that $\Gamma_{1,0} \land \Psi_{1,0} = (\Gamma_{1}
  \land \Psi_{1})_{0} = \Gamma_{0,0} \leq \Delta_{0,1}$.

  Next, we make use of Lemma \ref{lemma:useful} to prove that
  $\Gamma_{1,0} \land \Psi_{1,0} \not\leq \Delta_{0,1}$, thus
  obtaining a contradiction.
  By definition of the poset $\Cov(\Cov(L))$, $\Delta_{0} \lcover
  \Delta_{1}$, $\Gamma_{1} \leq \Delta_{1}$, $\Gamma_{1} \not\leq
  \Delta_{0}$ and, similarly, $\Psi_{1} \leq \Delta_{1}$ but $\Psi_{1}
  \not\leq \Delta_{0}$.
  In the statement of Lemma \ref{lemma:useful} let $\gamma =
  \Gamma_{1}$, $\delta = \Delta_{1}$, $\epsilon = \Delta_{0}$, and
  deduce $\Gamma_{1,0} \vee \Delta_{0,1} = \Delta_{1,1}$. Similarly,
  let $\gamma = \Psi_{1}$, $\delta = \Delta_{1}$, $\epsilon =
  \Delta_{0}$, and deduce $\Psi_{1,0} \vee \Delta_{0,1} =
  \Delta_{1,1}$. Therefore, from $\Gamma_{1,0} \vee \Delta_{0,1} =
  \Delta_{1,1} = \Psi_{1,0} \vee \Delta_{0,1}$ we deduce
  $(\Gamma_{1,0}
  \land \Psi_{1,0}) \vee \Delta_{0,1} = \Delta_{1,1}$ and, consequently,
  $\Gamma_{1,0} \land \Psi_{1,0} \not\leq \Delta_{0,1}$.
\end{proof}

\begin{proposition}
  \label{prop:CCcpb}
  If $L$ is join-semidistributive, then the projection
  $\pr[0]$ creates pullbacks.
\end{proposition}
\begin{proof}
  Since the projection $\pr[0]$ is a conservative Grothendieck
  fibration, it is enough by Lemma \ref{lemma:picreatespbs} to prove
  that if $\Gamma,\Psi\leq \Delta$, then there exists $\Upsilon \leq
  \Gamma,\Psi$ such that $\Upsilon_{0} = \Gamma_{0}\land\Psi_{0}$.

  To this goal, we observe first that $\Gamma_{0} \land \Psi_{0} <
  \Gamma_{1} \land \Psi_{1}$: as in proof the previous Proposition, we
  have $(\Gamma_{1} \land \Psi_{1})_{0} \vee \Delta_{0,1} =
  \Delta_{1,1}$, $(\Gamma_{1} \land \Psi_{1})_{0} \not\leq
  \Delta_{0,1}$, and consequently $\Gamma_{1} \land \Psi_{1} \not\leq
  \Delta_{0}$. Hence the standard relation $\Gamma_{0} \land \Psi_{0}
  \leq \Gamma_{1} \land \Psi_{1}$ is not an equality, since
  $\Gamma_{0} \land \Psi_{0} \leq \Delta_{0}$, and
  $\Gamma_{0} \land \Psi_{0} < \Gamma_{1} \land \Psi_{1}$.

  Let $(\Upsilon_{0},\Upsilon_{1}) \in \Cov^{2}(L)$ be such that $
  \Gamma_{0}\land\Psi_{0} = \Upsilon_{0} \lcover \Upsilon_{1} \leq
  \Gamma_{1} \land \Psi_{1}$ and observe that the following are
  equivalent: $\Upsilon_{1} \leq \Gamma_{0}$, $\Upsilon_{1} \leq
  \Delta_{0}$, $\Upsilon_{1} \leq \Psi_{0}$, $\Upsilon_{1} \leq
  \Gamma_{0} \land \Psi_{0} = \Upsilon_{0}$. For example, if
  $\Upsilon_{1} \leq \Gamma_{0}$, then $\Upsilon_{1} \leq \Gamma_{0}
  \leq \Delta_{0}$, and since $\Upsilon_{1} \leq \Psi_{1}$, then
  $\Upsilon_{1} \leq \Delta_{0} \land \Psi_{1} = \Psi_{0}$ and
  $\Upsilon_{1} \leq \Gamma_{0} \land \Upsilon_{0} = \Psi_{0}$.

  Since $\Upsilon_{1} \not\leq \Upsilon_{0}$, we deduce $\Upsilon_{1}
  \not\leq \Gamma_{0}$ and $\Upsilon_{1} \not\leq \Psi_{0}$.
  Therefore we have all the elements to state that $\Upsilon \leq
  \Gamma,\Delta$ and, by construction, we have $\Upsilon_{0} =
  \Gamma_{0} \land \Delta_{0}$.
\end{proof}

Using Proposition \ref{prop:CCcpb}, also in its dual form, we arrive
to the first achievement of this section.
\begin{theorem}
  \label{theo:derived}
  If $L$ is a semidistributive lattice, then $\Cov(L,\gamma)$ is a
  semidistributive lattice, for each $\gamma \in \Cov(L)$.
\end{theorem}

We shall use next the characterization of bounded lattices of Theorem
\ref{theo:bounded} to obtain the second main result of this section.
\begin{theorem}
  \label{theo:liftbounded}
  If $L$ is a bounded lattice, then so is $\Cov(L,\gamma)$, for each
  $\gamma \in \Cov(L)$.
\end{theorem}
\begin{proof}
  Since $L$ is bounded, then Theorem \ref{theo:bounded} ensures that
  it has a strict facet labelling $f : \Cov(L) \rTo \N$. Recall from
  Lemma \ref{lemma:lowercover} that a cover in $\Cov(L)$ is of the
  form $(\Gamma_{0},\Gamma_{1})$ with $\Gamma_{1} \pushdown[u]
  \Gamma_{0}$ for a unique lower cover $u$ of
  $\Gamma_{1,1}$. Therefore we define a function $F :
  \Cov(\Cov(L,\gamma)) \rTo \N$ by
  \begin{align*}
    F(\Gamma) & = f(u,\Gamma_{1,1})\,.
  \end{align*}
  Observe that if $\Gamma_{1} \pushdown[u] \Gamma_{0}$ and $\Gamma_{0}
  \pushup[w] \Gamma_{1}$, then $(u,\Gamma_{1,1})
  \spushdown[3pt]{\Gamma_{1,0}} (\Gamma_{0,0},w)$ so that
  $f(\Gamma_{0,0},w) = f(u,\Gamma_{1,1}) = F(\Gamma)$.  Therefore the
  definition of $F$ does not depend on whether we choose the relation
  $\pushdown$ or its dual $\pushup$.

  Next we prove that $F$ so defined is a strict facet
  labelling.  To this goal, let us suppose that $\Gamma
  \pushdown[\upsilon] \Delta$ and $\Gamma_{1} \leq \Upsilon_{0}
  \lcover \Upsilon_{1} \leq \upsilon$, as sketched in the next
  diagram:
  $$
  \xygraph{
    !~:{@{.}}
    []!c{{\Gamma_{\!1,}}}
    (
    [r(2)d]!c{{\Gamma_{\!0,}}}
    [d(3.1)l]!c{{\Delta_{\!0,}}}
    [l(1.2)u(1.1)]!c{{\Delta_{\!1,}}}
    ,
    [l(2)d(0.7)]!c{\upsilon}
    [d(0.6)r(0.4)]!c{{\Upsilon_{\!1,}}}
    [d(1.1)r(0.8)]!c{{\Upsilon_{\!0,}}}
    )
    !p{{\Gamma_{\!1,}}}{{\Gamma_{\!0,}}}{u}
    !e{{\Gamma_{\!1,}}}{\upsilon}
    !e{{\Gamma_{\!0,}}}{{\Delta_{\!0,}}}
    !u{{\Delta_{\!1,}}}{{\Delta_{\!0,}}}{z}
    !e{\upsilon}{{\Upsilon_{\!1,}}}
    !u{{\Upsilon_{\!1,}}}{{\Upsilon_{\!0,}}}{w}
    !e{{\Upsilon_{\!0,}}}{{\Delta_{\!1,}}}
  }
  $$
  Let us also suppose that $\Gamma_{1} \pushdown[u] \Gamma_{0}$, 
  $\Delta_{0} \pushup[z] \Delta_{1}$, and $\Upsilon_{0} \pushup[w]
  \Upsilon_{1}$, so that $F(\Gamma) = f(u,\Gamma_{1,1})$, $F(\Delta) =
  f(\Delta_{0,0},z)$, and $F(\Upsilon) = F(\Upsilon_{0,0},w)$.



  Recall that $\Delta_{0} = \upsilon \land \Gamma_{0}$ and therefore
  $\Delta_{0,0} = \upsilon_{0} \land \Gamma_{0,0} = \upsilon_{0} \land
  \Gamma_{1,0} \land u = \upsilon_{0} \land u$. 
  Since $\Delta_{0} \lcover z \leq \upsilon_{0}$, then $z \not\leq u$,
  otherwise $z \leq \upsilon_{0} \land u = \Delta_{0,0}$.  We
  have therefore $z \leq \Gamma_{1,1}$, $\Delta_{0,0} \leq u$, $z
  \not\leq u$, that is $(\Delta_{0,0},z) \leq (u,\Gamma_{1,1})$.
  Consequently, $F(\Delta) = f(\Delta_{0,0},z) = f(u,\Gamma_{1,1}) =
  F(\Gamma)$.
  

  In order to show that $F(\Delta) < F(\Upsilon)$ it is enough to show
  that if $j,k \in J(L)$, $(j_{\ast},j) \leq (\Upsilon_{0,0},w)$, and
  $(k_{\ast},k) \leq (u,\Gamma_{1,1})$, then $j A k$.
  It follows then, by Lemma  \ref{lemma:A}, that  
  \begin{align*}
    F(\Delta)  =   F(\Gamma) = f(u,\Gamma_{1,1}) & = f(k_{\ast},k) < f(j_{\ast},j)
    \leq f(\Upsilon_{0},w) = F(\Upsilon)\,.
  \end{align*}

  
  Let $m \in M(L)$ such that $(u,\Gamma_{1,1}) \leq (m,m^{\ast})$.  We
  have $j \leq w \leq \Gamma_{1,1} \leq m^{\ast}$, and if $j \leq
  m$, then $j \leq \upsilon_{0} \land \Gamma_{1,1} \land m =
  \upsilon_{0} \land u = \Delta_{0,0} \leq \Upsilon_{0,0}$, a
  contradiction.
  
  By duality, it also follows that $F(\upsilon,\Gamma_{1}) <
  F(\Upsilon)$.
\end{proof}


%% file: examples.tex
\section{Derived Lattices of Newman Lattices}
\label{sec:examples}

We refer the reader to \cite{bb,caspard,caspardbarbut,multinomial} for
introductory readings on Newman lattices.  In this section we
explicitly compute derived lattices $\Cov(L,\alpha)$ when $L$ is a
Permutohedron or an Associahedron and $\alpha \in \Cov(L)$ is an
\emph{atomic cover} of $L$, i.e. it is of the form $\alpha =
(\bot,\alpha_{1})$ with $\alpha_{1}$ an atom of $L$.
We shall see that these derived lattices are again Permutohedra
(respectively, Associahedra) of same dimension minus one. We remark
therefore a peculiar property of these lattices, they are
\emph{regular}, meaning that, up to isomorphism, the shape of
$\Cov(L,\alpha)$ does not depend on the choice of the atomic cover
$\alpha$.  We shall exhibit later a semidistributive lattice -- not
complemented -- that it not regular. 
Regularity is a reminiscent property of Boolean algebras: if ${\cal
  B}^{n}$ is the Boolean algebra with $n$ atoms, $\alpha_{1}$ being
one of them, then the equality $\Cov({\cal B}^{n},\alpha) = {\cal
  B}^{n-1}$ holds up to isomorphism. More generally:
\begin{proposition}
  If $\alpha$ is an atomic cover of a distributive lattice $L$, then
  the projection $\pr[0]$ from $\Cov(L,\alpha)$ to the lower set
  $\set{x \in L \mid \alpha_{1} \not\leq x}$ is an isomorphism.
\end{proposition}
The Proposition depends on modularity, since if $\alpha_{1} \not\leq
x$, then $x \lcover x \vee \alpha_{1}$.

\medskip

In the following proofs 
we shall intensively use the category of finite ordinals and functions
among them, a skeleton of the category of finite sets and
functions. To this goal, let $[n]$ be the set $\set{1,\ldots ,n}$ and,
for $i \in [n]$, denote by $\missi_{n} : [n-1] \rTo\relax [n]$ the
unique order preserving injection whose image is
$[n]\setminus\set{i}$.  For $k \in [n-1]$ denote by $\double{k}_{n} :
[n] \rTo \relax [n-1]$ the unique order preserving surjection such
that $\double{k}_{n}(k) = \double{k}_{n}(k + 1)$.  As the subscripts
$n$ will always be understood from the context, we shall omit them and
write only $\missi$ and $\double{k}$. 
%


\begin{proposition}
  \label{prop:derivperm}
  Let $\mathcal{S}_{n}$ be the Permutohedron on $n$ letters (i.e. the
  weak Bruhat order on permutations on $n$ elements).  If $\alpha$ is
  an atomic cover of $\mathcal{S}_{n}$ then
  $\Cov(\mathcal{S}_{n},\alpha)$ is isomorphic to $\mathcal{S}_{n
    -1}$.
\end{proposition}
\begin{proof}
  As usual we represent a permutation $w \in \mathcal{S}_{n}$ as the
  word $w(1)\ldots w(n) = w_{1}\ldots w_{n}$. An \emph{increase} of
  $w$ is an index $i \in \set{1,\ldots ,n-1}$ such that $w_{i} < w_{i
    + 1}$.  If $i$ is an increase of $w$ and $\sigma^{i}$ denotes the
  exchange permutation $(i,i+1)$, then we represent the cover $w
  \lcover w \circ \sigma^{i}$ of $\mathcal{S}_{n}$ by the pair
  $(w,i)$.  Every cover arises in this way.

  Remark next that a cover $(w,i)$ is perspective to the atomic cover
  $(\bot,\sigma^{k})$ if and only if $w_{i} = k$ and $w_{i + 1} = k +
  1$. If $(w,i)$ is such a cover, then we define $\psi^{k}(w,i)$ as
  the compose
  \begin{align*}
    &
    \mydiagram[4em]{%
      []*+{[n-1]}="S"
      :[u]*+{[n]}^{\missi}
      :[r(2.5)]*+{[n]}^{w}
      :[d]*+{[n-1]}="E"^{\double{k}}
      "S":"E"^{\psi^{k}(w,i)}
    }
  \end{align*}
  For example, $\psi^{2}(45231,3)$ is computed as follows.  We first
  erase the letter in third position and obtain the word
  $4531$. Then we normalize this word to a permutation. To this goal,
  knowing that in third position of the original word there was the
  letter $2$, we must decrease by one all the values of this word that
  are strictly greater than $2$. Thus, we obtain the permutation
  $3421$.
  We remark that we would have obtained the same result if we first
  erase the letter after the third position (i.e. in forth position)
  and then decrease by one all the values greater than $3$. More
  generally, we could have equivalently defined $\psi^{k}(w,i)$ as the
  compose $\double{k + 1} \circ w \circ \miss{i + 1}$.

  The informal example already suggests that $\psi^{k}(w,i)$ is
  injective, and hence it is bijective; let us argue formally in this
  sense. If $\psi^{k}(w,i)$ is not injective, then there exists $x, y
  \in [n]$ such that $x, y,i$ are pairwise distinct and
  $\double{k}(w_{x})= \double{k}(w_{y})$.  But this may happen only if
  $\set{w_{x},w_{y}} = \set{k,k+1}$ and, by the assumption on the
  cover $(w,i)$ stating that $w_{i} = k$ and $w_{i + 1} = k +1$, this
  happens exactly when $\set{x,y} = \set{i,i+1}$.
  

  \medskip It is easily seen that $\psi^{k}$ is a bijection from
  $\Cov({\cal S}_{n}, (\bot,\sigma^{k}))$ to ${\cal S}_{n-1}$: if $u
  \in {\cal S}_{n-1}$, then there exists a unique cover $(w,i)$ of
  ${\cal S}_{n}$ which is sent by $\psi^{k}$ to $u$. The position $i
  \in [n-1]$ is the unique index $i$ such that $u_{i} = k$ and then we
  define $w$ as follows:
  \begin{align*}
    w_{j} & =
    \begin{cases}
      k  \,,& j = i \,,\\
      \miss{k }(u_{\double{i}(j)})\,,
      &\toth\,.
    \end{cases}
  \end{align*}


  To prove that $\psi^{k}$ is an order isomorphism, we prove that
  $(w,i) \lcover (w',i')$ if and only if $\psi^{k}(w,i) \lcover
  \psi^{k}(w',i')$.  This equivalence is an immediate consequence of
  the following two claims.
  \begin{claim}
    For $j \in [n-2]$, $j$ is an increase of $\psi^{k}(w,i)$ if and
    only if $\missi(j)$ is an increase of $w$.
  \end{claim}
  Since $w_{i} = k$, then $w_{\missi(j)} \neq k$ for each $j \in
  [n-1]$, so that $\miss{k}$ and $\double{k}$ are inverse order
  preserving bijections relating the sets $[n-1]$ and
  $\set{w_{\missi(j)} \mid j \in [n-1]} = [n] \setminus \set{k}$.
  Thus we have
  \begin{align*}
    \psi^{k}(w,i)_{j} < \psi^{k}(w,i)_{j + 1}
    & \tiff
    \double{k}(w_{\missi(j)}) < \double{k}(w_{\missi(j + 1)})
    \\
    & \tiff
    w_{\missi(j)} <  w_{\missi(j +1)}\,,
  \end{align*}
  and we are left to argue that
  \begin{align*}
    w_{\missi(j)} <  w_{\missi(j +1 )}
    & \tiff w_{\missi(j)} < w_{\missi(j) + 1}\,.
  \end{align*}
  Clearly this is the case if $\missi(j) + 1 = \missi(j + 1)$.
  If not, then $\missi(j) = i -1$, and we need to argue that $w_{i
    -1} < w_{i}$ if and only if $w_{i-1} < w_{i + 1}$. This is an
  immediate consequence of the assumptions $w$, namely that $w_{i} =
  k < k+1 = w_{i + 1}$.
  \eproofofclaim

  

  The previous claim implies that the following data are in bijection:
  the upper covers of $\psi^{k}(w,i)$, the increases of
  $\psi^{k}(w,i)$, the increases of $w$ that are distinct of $i$, the
  upper covers of $w$ that are distinct from $w \circ \sigma^{i}$, the
  upper covers of $(w,i)$ in the poset $\Cov({\cal S}_{n})$. The next
  claim shows that $\psi^{k}$ respects this bijection.
  \begin{claim}
    Given an increase $j$ of $\psi^{k}(w,i)$, let $w'$ and $i'$ be
    determined by the pushup relation $(w,i) \pushup[w \circ
    \sigma^{\missi(j)}] (w',i')$. Then
    \begin{align*}
      \psi^{k}(w',i') & = \psi^{k}(w,i) \circ \sigma^{j}\,.
    \end{align*}
  \end{claim}
  Firstly, we claim first that the following relations hold: 
  
  \begin{align}
    \missi \circ \sigma^{j}
    & =   
    \sigma^{\missi(j)} \circ \missi
    \,,
    && \!\tif |\missi(j) -i | > 1  \,,
    \label{eq:ij}
    \\
    \missi \circ \sigma^{j} & = \sigma^{\missi(j)} \circ
    \sigma^{i} \circ \miss{\missi(j)}\,, && \toth, \tif
    |\missi(j) -i | =1\,.
    \label{eq:ijij}
    \end{align}
    We leave the reader to verify equation \eqref{eq:ij} and focus
    instead on \eqref{eq:ijij}. We shall split its verification into
    two obvious cases and suggest a proof by making use of standard
    strings diagrams, see for example \cite[\S 2.3]{lafont}.


    Case $\missi(j) = i + 1$, that is $j = i$. Equation
    \eqref{eq:ijij} reduces to $\missi \circ \sigma^{i} = \sigma^{i
      + 1} \circ \sigma^{i} \circ \miss{i + 1}$, which is easily seen
    to hold by comparing  the following two diagrams:
    $$
    \mydiagram[3em]{%
      []!h
       "i0"-@`{"j1"}"k2"
       "i1"[l(0.5)]="ii1"
       "j0"-@`{"ii1"}"j2"
       "i1"!k-"i2"
    }
    \;\;=\;\;
    \mydiagram[3em]{%
      []!g
      "i0"-@`{"i1","j2"}"k3"
      "j0"-@`{"k1","k2"}"j3"
      "j0"!k-@`{"j1","i2"}"i3"
    }
    $$
    
    Case $\missi(j) = i - 1$, that is $j = i -1$.  Equation
    \eqref{eq:ijij} reduces to $\missi \circ \sigma^{i -1} =
    \sigma^{i - 1} \circ \sigma^{i} \circ \miss{i - 1}$, which is
    easily seen to hold by comparing  the following two diagrams:
    $$
    \mydiagram[3em]{%
      []!H
       "i0"-@`{"j1"}"k2"
       "j0"-@`{"i1"}"i2"
       "j1"!k-"j2"
    }
    \;\;=\;\;
    \mydiagram[3em]{%
      []!G
      "i0"!k-@`{"i1","i2"}"j3"
      "i0"-@`{"j1","k2"}"k3"
      "j0"-@`{"k1","j2"}"i3"
    }
    $$
    
    \medskip 

    Next, we compute $w'$ and $i'$ in the pushup $(w,i)
    \pushup[w \circ \sigma^{\missi(j)}] (w',i')$.
    If $|\missi(j) - i| > 1$, then $w' = w \circ
    \sigma^{\missi(j)}$ and $i' = i$.
    Otherwise, if $|\missi(j) - i| = 1$, then $w ' = w \circ
    \sigma^{\missi(j)} \circ \sigma^{i}$ and $i' = \missi(j)$.

    \medskip

    Therefore, if $|\missi(j) - i| > 1$, then
    \begin{align*}
      \psi^{k}(w',i')
      & = \psi^{k}(w \circ \sigma^{\missi(j)},i)
      \\
      & = \double{k}\circ w \circ \sigma^{\missi(j)} \circ \missi
      \\
      & = \double{k}\circ w \circ \missi \circ \sigma^{j}
      = \psi^{k}(w,i) \circ \sigma^{j}\,,
      \intertext{and,  if $|\missi(j) - i| = 1$, then}
      \psi^{k}(w',i')
      & = \psi^{k}(w \circ \sigma^{\missi(j)}\circ \sigma^{i},\missi(j))
      \\
      & = \double{k}\circ w \circ \sigma^{\missi(j)} \circ
      \sigma^{i} \circ \miss{\missi(j)}
      \\
      & = \double{k}\circ w \circ \missi \circ \sigma^{j}
      = \psi^{k}(w,i) \circ \sigma^{j}\,.
      \tag*{\eproofofclaimlabel}
     \end{align*}
  This also completes the proof of Proposition \ref{prop:derivperm}. 
\end{proof}

\medskip

We use Proposition \ref{prop:derivperm} to argue that derived
semidistributive lattices of the form $\Cov(L,\gamma)$ are not
quotients of $L$ in the most obvious way.  It is a standard reasoning
to argue that $\delta \in \Cov(L,\gamma)$ implies
$(\delta_{0},\delta_{1}) \in \theta(\gamma_{0},\gamma_{1})$, where
$\theta(\gamma_{0},\gamma_{1})$ is the congruence generated by the
pair $(\gamma_{0},\gamma_{1})$. It is reasonable to ask whether the
lattice $\Cov(L,\gamma)$ is related to the specific quotient lattice
$L/\theta(\gamma_{0},\gamma_{1})$.  The following Proposition gives a
first answer in the negative, showing that these two lattices are not
in general isomorphic.
\begin{Proposition}
  For $k \in \set{1,\ldots ,n-1}$ the lattice
  $\mathcal{S}_{n}/\theta(\bot,\sigma^{k})$ is isomorphic to the
  lattice $\mathcal{S}_{k}\times \mathcal{S}_{n-k}$.
\end{Proposition}
\begin{proof}
  We recall that for $L$ a lattice $L$ and $\theta$ a congruence  of
  $L$, each equivalence class $[x]_{\theta}$ has a least element
  $\mu_{\theta}(x)$, computed as follows:
  \begin{align*}
    \mu_{\theta}(x) & = \bigvee \set{j \in J(L)\mid j \leq x \tand
      (j_{\ast},j) \not\in \theta }\,.
  \end{align*}
  The quotient $L/\theta$ is then isomorphic to the poset $\langle
  \set{\mu_{\theta}(x) \mid x \in L},\leq\rangle$. We use this
  representation to give explicit form to
  $\mathcal{S}_{n}/\theta(\bot,\sigma^{k})$.

  
  Recall that a permutation is join-irreducible iff it has a unique
  descent, i.e. a unique index $i \in \set{1,\ldots ,n-1}$ such that
  $w_{i} > w_{i + 1}$.  In \cite{multinomial} we called the pair
  $(w_{i + 1},w_{i})$ the principal plan of the join-irreducible $w$.
  Let us denote by $\Dclosure$ the reflexive transitive closure of the
  join-dependency relation between join-irreducible permutations.
  From the characterization given there of the join-dependency
  relation, we have that $w \Dclosure \sigma^{k}$ iff the $[k,k+1]
  \subseteq [a,b]$, where $a,b$ is the principal plan of $w$.
  Therefore, for $w$ a join-irreducible permutation, we have the
  following equivalences: $(w_{\ast},w) \in \theta(\bot,\sigma^{k})$,
  iff $w \Dclosure \sigma^{k}$ iff iff $[k,k+1] \subseteq [a,b]$,
  $a,b$ the principal plan of $w$.
  Consequently, a join-irreducible permutation is not congruent to its
  unique lower cover modulo $\theta(\bot,\sigma^{k})$ if and only if
  its principal plan does not contain the interval $[k,k+1]$.

  Remark that, for $a < b < c$, we have $[k,k+1] \subseteq [a,c]$ if
  and only if $[k,k+1] \subseteq [a,b]$ or $[k,k+1] \subseteq [b,c]$.
  Using this fact, we see that if $D = D(w)$ is the set of
  disagreements (or inversions) of some permutation $w$, then $D' =
  \set{(a,b) \in D \mid [k,k+1] \not\subseteq [a,b] }$ is also the set
  of disagreements of some permutation $w'$.  To this goal, it is
  enough to verify that $D'$ is closed -- i.e. $(a,b),(b,c) \in D'$
  implies $(a,c) \in D'$ -- and open as well -- i.e. $a < b < c$ and
  $(a,c) \in D'$ implies $(a,b) \in D'$ or $(b,c) \in D'$.  Since
  $(a,b) \in D(w)$ if and only if there exists $j \in J({\cal S}_{n})$
  such that $(a,b)$ is the principal plan of $j$, then we deduce that
  the permutation $w'$ is the least element in the congruence class of
  $w$.
  
  Knowing that the order on ${\cal S}_{n}$ is given by inclusion of
  disagreement sets, the relation 
  \begin{align*}
    D(w') & = \set{(a,b) \in D(w) \mid [a,b] \subseteq [1,k]} \uplus 
    \set{(a,b) \in D(w) \mid [a,b] \subseteq [k + 1,n]}
  \end{align*}
  exhibits ${\cal S}_{n}/\theta(\bot,\sigma^{k})$ as the order
  theoretic product of ${\cal S}_{k}$ and ${\cal S}_{n-k}$.
\end{proof}
Considering that finite semidistributive lattices and bounded lattices
form pseudovarieties \cite{jonsonnrival,nation}, we leave it as an
open problem whether derived lattices are constructible by means of
standard operations such as homomorphic images, subalgebras and
products.

\medskip

We are ready to tackle computation of the lattices derived from
Associahedra by atomic covers. The computation we present here is a
direct one. Considering however that Associahedra are quotient of
Permutohedra, see \cite[\S 9]{reading2}, we expect that the next
Proposition may be derived from Proposition \ref{prop:derivperm} in a
more informative manner.
\begin{proposition}
  \label{prop:derivasso}
  Let $\mathcal{T}_{n}$ be the Associahedron on $n+1$ letters (i.e.
  the Tamari lattice).  If $\alpha$ is an atomic cover of
  $\mathcal{T}_{n}$ then $\Cov(\mathcal{T}_{n},\alpha)$ is isomorphic
  to $\mathcal{T}_{n -1}$.
\end{proposition}
To prove Proposition \ref{prop:derivasso}, we shall review some facts
about the explicit representation of the Tamari lattices as lattices
of bracketing vectors with the pointwise order, see
\cite{huangtamari,bb,caspardbarbut}.  A \emph{bracketing vector} is a
vector $v \in \set{1,\ldots ,n}^{n}$ such that (i) $i \leq v_{i}$ and
(ii) $i < j \leq v_{i}$ implies $v_{j} \leq v_{i}$.  We are going to
determine covers of the pointwise order.  Let us say that $k \in
\set{1,\ldots ,n-1}$ is a \emph{split} of a bracketing vector $v$ if
$i < k \leq v_{i}$ implies $v_{v_{k} + 1} \leq v_{i}$. The next
diagram should help understanding the condition: if $k$ is a split,
then the black region is a forbidden area meaning that it does not
contain points of the form $(i,v_{i})$.
$$
\mydiagram[3em]{
  []!t{1}{5}
  "SE1"[r(1.2)u(1.2)]!t{k}{1.5}
  "SE1"[r(3)u(3)]!t{{v_{k}+1}}{1}
  "SE1"[u(2.7)l(0.5)]*+{v_{k}}
  "SE1"[u(4)]="A"[l(0.5)]*+{v_{v_{k}+1}}
  "NOk"[u(0.35)l(0.35)]="B"
  "B"[l(0.1)]="B"
  !{{"A"."B"}*[mac28]++\frm{*}}
}
$$

\begin{lemma}
  Let $k$ be a split of a bracketing vector $v$ and define the vector
  $v^{k}$ by
  \begin{align*}
    v^{k}_{i} & =
    \begin{cases}
      v_{v_{k} + 1}, & i = k,\\
      v_{i}, & \text{otherwise}.
    \end{cases}
  \end{align*}
  Then $v \lcover v^{k}$ and moreover all the covers in ${\cal T}_{n}$
  arise in this way.
\end{lemma}
\begin{proof}
  We observe first that $v^{k}$ is again a bracketing vector.
  Condition (i) is satisfied: if $i \neq k$, then $i \leq v_{i} =
  v^{k}_{i}$, and otherwise $k \leq v_{k} < v_{k} + 1 \leq v_{v_{k} +
    1} = v^{k}_{k}$.
  Condition (ii) clearly holds if both $i$ and $j$ are distinct from
  $k$. Let us suppose that $j =k$, that is $i < k \leq v^{k}_{i} =
  v_{i}$. Since $k$ is a split, then $v^{k}_{k} = v_{v_{k} +1} \leq
  v_{i}$.  Let us suppose that $i = k$, that is, $k < j \leq
  v^{k}_{k} = v_{v_{k} +1}$.  If $j \leq v_{k}$ then $v_{j} \leq v_{k}
  < v^{k}_{k}$.  If $v_{k} < j$, then $v_{k} + 1 \leq j \leq v_{v_{k}
    +1}$ and $v_{j} \leq v_{v_{k} + 1} = v^{k}_{k}$.

  Let us suppose that $v < w$ and let $k$ be the least index such that
  $v_{k} < w_{k}$. Observe first that $v_{v_{k} + 1} \leq w_{k}$: from
  $v_{k} < w_{k}$ we can write $k < v_{k} + 1 \leq w_{k}$ and hence
  $v_{v_{k} + 1} \leq w_{v_{k} + 1} \leq w_{k}$. Also $k$ is a split
  of $v$: if $i < k \leq v_{i}$ then $i < k \leq w_{i}$, $w_{k} \leq
  w_{i}$ so that $v_{v_{k} + 1} \leq w_{k} \leq w_{i } = v_{i}$. This
  shows that $v \lcover v^{k}$ and moreover that any upper cover of
  $v$ is of the form $v^{k}$ for some split of $v$.
\end{proof}

\begin{lemma}
  \label{lemma:Tnfactes}
  Let $v$ be a bracketing vector.
  \begin{enumerate}
  \item If $i,j$ are two splits of $v$ with $j = v_{i} +
    1$, then we have the following pentagon:
    $$
    \xygraph{%
      []*+{v_{i}\ldots v_{j} \ldots v_{v_{j} +1}}
      (-[ul]*+{v_{j}\ldots v_{j} \ldots v_{v_{j} +1}}|{i} 
      -[u]*+{v_{v_{j} + 1}\ldots v_{j} \ldots v_{v_{j} +1}}|{i} 
      -[ur]*+{v_{v_{j} + 1}\ldots v_{v_{j} +1} \ldots v_{v_{j}
          +1}}="E"|{j} ) 
      -[u(1.5)r]*+{v_{i}\ldots v_{v_{j} +1} \ldots v_{v_{j} +1}}|{j}
      -"E" |{i}}
    $$
  \item If $i,j$ are two splits of $v$ and $j \neq v_{i} + 1$, then we
    have one of the following diamonds:
    $$
    \xygraph{%
      []*+{v_{i}\ldots v_{v_{i} +1} \ldots v_{j} \ldots v_{v_{j} +1}}
      (-[ull]*+{v_{v_{i} +1}\ldots v_{v_{i} +1} \ldots v_{j} \ldots
        v_{v_{j} +1} } |{i} 
      -[urr]*+{v_{v_{i} +1}\ldots v_{v_{i} +1} \ldots v_{v_{j} +1} \ldots
        v_{v_{j} +1} }="E"|{j}) 
      -[urr]*+{v_{i}\ldots v_{v_{i} +1} \ldots v_{v_{j} +1}
        \ldots v_{v_{j} +1}}|{j} 
      -"E" |{i}}
    $$
    $$
    \xygraph{%
      []*+{v_{i}\ldots v_{j} \ldots v_{v_{j} +1} \ldots  v_{v_{i} +1} }
      (-[ull]*+{v_{v_{i} +1}\ldots v_{j} \ldots
        v_{v_{j} +1} \ldots v_{v_{i} +1}  }|{i} 
      -[urr]*+{v_{v_{i} +1}\ldots v_{v_{j} +1} \ldots
        v_{v_{j} +1} \ldots v_{v_{i} +1}  }="E"|{j}) 
      -[urr]*+{v_{i}\ldots v_{v_{j} +1}
        \ldots v_{v_{j} +1} \ldots v_{v_{i} +1} }|{j} 
      -"E"|{i} }
    $$
  \end{enumerate}
\end{lemma}
We refer the reader to \cite[Propositions 4 and 5]{caspardbarbut} for
a detailed proof of this Lemma.

If $v$ is a bracketing vector and $k$ is a split of $v$, then we
denote the cover $v \lcover v^{k}$ by the pair $(v,k)$. From the
previous Lemma it immediately follows:
\begin{corollary}
  If $i \neq j$ and $(v,j) \pushup[v^{i}] (w,k)$, then $k = j$ and
  \begin{align*}
    w & =
    \begin{cases}
      v^{ii}, & j = v_{i} + 1, \\
      v^{i}, & \text{otherwise}.
    \end{cases}
  \end{align*}
\end{corollary}
\begin{corollary}
  \label{cor:unique}
  A cover $(v,j)$ is perspective to the atom $(\bot,k)$ if and only if
  $j = k$ and $v_{k} = k$.  In this case $k$ is the unique index $i$
  such that $v_{i} = k$.
\end{corollary}
\begin{proof}
  The condition is necessary: the property holds for $(v,k)$ and, by
  Lemma \ref{lemma:Tnfactes}, is preserved under the operation of
  pushing up covers. The condition is also sufficient, for which it is
  enough to remark that if $v_{k} = k$, then $\bot \leq v$, $\bot^{k}
  = k + 1\not\leq k = v_{k}$, $\bot^{k} = k + 1 \leq v^{k}_{k} =
  v_{v_{k} + 1}$.
  For the last statement, let us suppose that $k$ is a split of $v$
  and that $v_{k} = k$. If $v_{i} = k$, then $i \leq k$. However $i <
  k$ contradicts $i$ being a split.
\end{proof}

We are ready to proof Proposition \ref{prop:derivasso}. 
  
\begin{proof}[Proof of Proposition \ref{prop:derivasso}]
  From a cover $(v,k) \in \Cov(\mathcal{T}_{n})$, perspective to an
  atom, define a bracketing vector $\psi(v,k) \in \mathcal{T}_{n-1}$
  as the compose
  $$
  \xygraph{%
    []*+{[n-1]}="S"
    :[u]*+{[n]}^{\miss{k}}
    :[r(2.5)]*+{[n]}^{v}
    :[d]*+{[n-1]}="E"^{\double{k}}
    "S":"E"^{\psi(v,k)}
    }\,.
  $$
  Before carrying on with the proof, we collect first some
  remarks. Observe that $\double{k}(\miss{k}(x)) = x$, while $x \leq
  \miss{k}(\double{k}(x))$ and this is an equality if $x \neq k$.
  Therefore $\miss{k}$ is right adjoint to $\double{k}$ and moreover
  $\double{k}$ is inverse to $\miss{k}$ if restricted to $[n]\setminus
  \set{k}$. Also  $\double{k}(x) + 1 = \double{k}(x + 1)$
  if $x \neq k$.  If $(v,k)$ is perspective to an atom, so that $v_{j}
  = k$ implies $j = k$, an integer of the form $v_{\miss{k}(j)}$ is not
  equal to $k$, otherwise $k = \miss{k}(j)$, a contradiction.
  Consequently we shall use formulas such as
  $\miss{k}(\double{k}(v_{\miss{k}(j)})) = v_{\miss{k}(j)}$, and
  $\double{k}(v_{\miss{k}(j)}) + 1 = \double{k}(v_{\miss{k}(j)} + 1)$.


  Let us verify that $\psi(v,k)$ is a bracketing vector. The relation
  $i \leq \psi(v,k)_{i} = \double{k}(v_{\miss{k}(i)})$ immediately
  follows from $\miss{k}(i) \leq v_{\miss{k}(i)}$.  If $i < j \leq
  \psi(v,k)_{i} = \double{k}(v_{\miss{k}(i)})$, then $\miss{k}(i) <
  \miss{k}(j) \leq \miss{k}(\double{k}(v_{\miss{k}(i)})) =
  v_{\miss{k}(i)}$ and $v_{\miss{k}(j)} \leq v_{\miss{k}(i)}$ since
  $v$ is a bracketing vector; the relation $\psi(v,k)_{j} \leq
  \psi(v,k)_{i}$ follows then by applying $\double{k}$.

  The correspondence $\psi$ is a bijection: given $w \in
  \mathcal{T}_{n-1}$ the vector $v \in \mathcal{T}_{n}$, defined by
  $v_{i} = k$ if $i = k$ and $v_{i} = \miss{k}(w_{\double{k}(i)})$ otherwise,
  is the unique bracketing vector such that $(v,k)$ is a cover
  perspective to $(\bot,k)$ and $\psi(v,k) = w$.
  
  We are going to verify that (a) $j$ is a split of $\psi(v,k)$ iff
  $\miss{k}(j)$ is a split of $v$, (b) if $(v,k) \pushup[v^{j}] (w,k)$
   then $\psi(w,k) =
  \psi(v,k)^{\double{k}(j)}$.  From these properties it follows that
  $\psi$ preserves and reflects the covering relation and therefore
  it is an order isomorphism.

  (a) Let us suppose first that $\miss{k}(j)$ is a split of $v$ and
  that $l < j \leq \psi(v,i)_{l} = \double{k}(v_{\miss{k}(l)})$. It
  follows that $\miss{k}(l) < \miss{k}(j) \leq \miss{k}(\double{k}(
  v_{\miss{k}(l)})) = v_{\miss{k}(l)}$ and therefore
  $v_{v_{\miss{k}(j)} + 1} \leq v_{\miss{k}(l)}$. Hence
  \begin{align*}
    \psi(v,k)_{\psi(v,k)_{j} + 1} & =
    \double{k}(v_{\miss{k}(\double{k}(v_{\miss{k}(j)}) + 1)}) =
    \double{k}(v_{\miss{k}(\double{k}(v_{\miss{k}(j)} + 1))})
    \\
    & = \double{k}(v_{v_{\miss{k}(j)} +1}) \leq
    \double{k}(v_{\miss{k}(l)}) = \psi(v,k)_{l}\,.
  \end{align*}
  Let us suppose now that $j$ is a split of $\psi(v,k)$ and that $l <
  \miss{k}(j) \leq v_{l}$. Observe that the relation $l < v_{l}$
  implies that $l \neq k$. Since both $l$ and $\miss{k}(j)$ are
  distinct from $k$, the relation $l < \miss{k}(l)$ is strictly
  preserved by $\double{k}$ and consequently
  \begin{align*}
    \double{k}(l) & < \double{k}(\miss{k}(j)) = j \\
    & \leq \double{k}(v_{l}) =
    \double{k}(v_{\miss{k}(\double{k}(l))}) =
    \psi(v,k)_{\double{k}(l)}\,.
  \end{align*}
  We have therefore $\psi(v,k)_{\psi(v,k)_{j} + 1} \leq
  \psi(v,k)_{\double{k}(l)}$ and
  \begin{align*}
    \double{k}(v_{v_{\miss{k}(j)} + 1})
    & =
    \double{k}(v_{\miss{k}(\double{k}(v_{\miss{k}(j)}    + 1 ))}) 
    =
    \double{k}(v_{\miss{k}(\double{k}(v_{\miss{k}(j)} )   + 1)}) \\
    & =
    \psi(v,k)_{\psi(v,k)_{j} + 1}
    \leq \psi(v,k)_{\double{k}(l)}
    = \double{k}(v_{l})\,.
  \end{align*}
  Transposing this relation and considering that $l \neq k$ we deduce
  \begin{align*}
    v_{v_{\miss{k}(j)} + 1} 
    &\leq \miss{k}(\double{k}(v_{l})) = v_{l}\,.
  \end{align*}
  
  (b) Let us suppose that $(v,k) \pushup[v^{j}] (w,k)$, so that $w =
  v^{jj}$ if $k = v_{j} + 1$ and $w = v^{j}$ otherwise.  We
  want to prove that $\psi(w,k) = \psi(v,k)^{\double{k}(j)}$. Let us
  begin to show that these two vectors coincide in each component $i$
  such that $i \neq \double{k}(j)$ (or equivalently $j \neq
  \miss{k}(i)$):
  \begin{align*}
    \psi(w,k)_{i} & = \double{k}(w_{\miss{k}(i)}) =
    \double{k}(v_{\miss{k}(i)})
    = \psi(v,k)_{i} = \psi(v,k)_{i}^{\double{k}(j)}\,.
  \end{align*}
  Therefore we are left to compare the values of the two vectors at
  the coordinate $i = \double{k}(j)$. On the one hand, we have
  \begin{align*}
    \psi(v,k)_{\double{k}(j)}^{\double{k}(j)} & =
    \psi(v,k)_{i}^{i}
    = \psi(v,k)_{\psi(v,k)_{i}  + 1}
    = \double{k}(v_{\miss{k}(\double{k}(v_{\miss{k}(i)}) + 1)})
    \\&
    = \double{k}(v_{\miss{k}(\double{k}(v_{\miss{k}(i)} + 1))})
    = \double{k}(v_{\miss{k}(\double{k}(v_{j} + 1))})
    \\&
    =
    \begin{cases}
      \double{k}(v_{k +1})\,, & k = v_{j} + 1 \\
      \double{k}(v_{v_{j} + 1})\,, & \toth.
    \end{cases}
  \end{align*}
  On the other hand, we have
  \begin{align*}
    \psi(w,k)_{\double{k}(j)} & 
    = \\
    \double{k}(w_{j}) 
    &
    =
    \begin{cases}
      \begin{array}[b]{@{\hspace{0mm}}l@{\hspace{0mm}}l}
        \double{k}(v^{jj}_{j})
        & = \double{k}(v^{j}_{v^{j}_{j} + 1})
        = \double{k}(v^{j}_{v_{v_{j} + 1} + 1}) 
        \\
        &
        = \double{k}(v^{j}_{v_{k} + 1})
        =\double{k}(v^{j}_{k + 1})
        \\
        & =\double{k}(v_{k + 1})\,,
      \end{array} 
       & k = v_{j} + 1\,, \\[4mm]
       \double{k}(v^{j}_{j}) = \double{k}(v_{v_{j} + 1})\,, & \toth .
     \end{cases}
   \end{align*}
   This completes the proof of Proposition \ref{prop:derivasso}.
 \end{proof}

 Let us say that a finite semidistributive lattice is \emph{regular}
 if the lattices $\Cov(L,\alpha)$, $\alpha$ an atomic cover of $L$,
 are all isomorphic.  It is not the case that every semidistributive
 lattice is regular as witnessed for example by the multinomial
 lattice $\mathcal{L}(2,2,1)$. The reader can find in \cite[\S 6]{bb}
 an introductory discussion of multinomial lattices. Let us recall
 that, for $v \in \N^{k}$, the elements of the multinomial lattice
 ${\cal L}(v)$ are words $w$ over an ordered alphabet
 $\set{a_{1},\ldots ,a_{k}}$ such that, for $i = 1,\ldots ,k$, the
 number $|w|_{a_{i}}$ of occurrences of the letter $a_{i}$ equals
 $v_{i}$. The
 Hasse diagram of ${\cal L}(v)$ is obtained by exchanging the position
 of contiguous letters that appear in the right order.
 The bottom of the lattice $\mathcal{L}(2,2,1)$ is represented in
 figure 1. 
 Let $\alpha = aaabbc \lcover ababc$ and $\beta = aaabbc \lcover
 aabcb$ be two atoms of this lattice, if we consider the bottoms of
 $\Cov(L,\alpha)$ and $\Cov(L,\beta)$ we observe these two lattices
 are not isomorphic. We remark that the lattice $\mathcal{L}(2,2,1)$
 is not complemented, contrary to the Newman lattices considered in
 this section. It might be conjectured that complemented
 semidistributive lattices are regular. More generally it is an open
 problem to identify sufficient conditions that ensure that a
 semidistributive lattice is regular.
\begin{figure}[h]
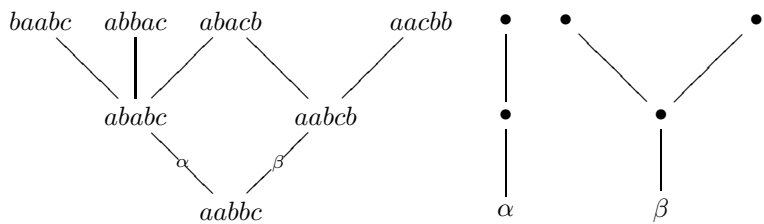

  \label{fig:multinomial}
  $$
  \xygraph{ []*+{aabbc}="aabbc" ([lu]*+{ababc}="ababc"
    ([lu]*+{baabc}="baabc",[u]*+{abbac},[ru]*+{abacb}),
    [ru]*+{aabcb}="aabcb" [ru]*+{aacbb}="aacbb" )
    "aabbc"(-"ababc"|{\alpha},-"aabcb"|{\beta})
    "ababc"(-"baabc",-"abbac",-"abacb") "aabcb"(-"abacb",-"aacbb") }
  \;\;\;\;
  \xygraph{
    []*+{\alpha}-[u]*+{\bullet}-[u]*+{\bullet}
  }
  \;\;\;\;
  \xygraph{
    []*+{\beta}-[u]*+{\bullet}(-[lu]*+{\bullet},-[ru]*+{\bullet})
  }
  $$\centering
  \caption{The bottom of the lattices $L = \mathcal{L}(2,2,1)$,
    $\Cov(L,\alpha)$ and $\Cov(L,\beta)$.}
\end{figure}
